\documentclass[10pt, oneside]{article}
\usepackage{mathrsfs}
\usepackage{amsfonts}
\usepackage{amssymb}
\usepackage{amsmath,amsfonts,amssymb,amsthm}
\usepackage{mathtools}
\usepackage{caption}

\usepackage{cite}

\usepackage{tikz}
\usepackage{subfigure}
\usepackage{float}
\usepackage{paralist}
\usepackage{indentfirst}
\usepackage{authblk}
\usepackage[colorlinks,
linkcolor=red,
citecolor=blue
]{hyperref}
\usepackage{color}
\usepackage{graphicx}
\usepackage{extarrows}

\usepackage{dsfont}

\usepackage[utf8]{inputenc} 
\usepackage[T1]{fontenc}
\usepackage{lmodern}       
\input{glyphtounicode}

\DeclareMathOperator{\dive}{div}

\numberwithin{equation}{section}

\def\bma#1\ema{{\allowdisplaybreaks\begin{aligned}#1\end{aligned}}}

\numberwithin{equation}{section}

\usepackage[title,titletoc]{appendix}

\usepackage{titlesec}

\newcommand{\B}{{\mathbb B}}

\newcommand{\N}{{\mathbb N}}

\newcommand{\R}{{\mathbb R}}

\newcommand{\Z}{{\mathbb Z}}

\newcommand{\sC}{{\mathcal C}}

\newcommand{\sL}{{\mathcal L}}

\newcommand{\sS}{{\mathcal S}}

\newtheorem{prop}{Proposition}[section]

\def\dz1{\nabla Z_1}
\def\dive{ \hbox{\rm div}\,  }

\definecolor{red}{rgb}{1.00,0.00,0.00}
\definecolor{blue}{rgb}{0.00,0.00,0.63}
\definecolor{black}{rgb}{0.00,0.00,0.00}
\definecolor{purple}{rgb}{0.00,1.00,0.00}
\definecolor{pink}{rgb}{0.95,0.01,0.08}

\usepackage{xcolor}

\newtheorem{theorem}{Theorem}[section]

\newtheorem{lemma}{Lemma}[section]

\newtheorem{remark}{Remark}[section]

\allowdisplaybreaks
\def\dive{ \hbox{\rm div}\,  }
\def\ddj{\dot{\Delta}_{j}}

\def\var{\varepsilon}

\def\dZ_1{\delta\!Z_1}

\def\div{ \hbox{\rm div}\,  }

\numberwithin{equation}{section}

\begin{document}

\title{Compressible Euler equations with time-dependent damping in the critical regularity setting: global well-posedness and strong relaxation limit}

\author{Timothée Crin-Barat, Xinghong Pan, Ling-Yun Shou \& Qimeng Zhu}

\date{}

\maketitle
\begin{abstract}
We investigate the relaxation problem and the diffusion phenomenon for the compressible Euler system with time-dependent damping coefficients of the form $\tfrac{\mu}{(1+t)^{\lambda}}$ in $\mathbb{R}^d$ $(d \geq 1)$. We establish uniform regularity estimates with respect to the relaxation parameter $\var$ and prove the global well-posedness of classical solutions to the Cauchy problem for initial data around equilibrium in critical hybrid Besov spaces. 
In addition, we justify the global-in-time strong convergence of the solutions toward those of a general porous medium-type diffusion system, for ill-prepared initial data, with an explicit convergence rate. The core of our proof relies on a refined hypocoercivity framework combined with a new time-dependent frequency decomposition, both adapted to handle the non-autonomous damping structure. This enables us to treat the overdamped regime $\lambda \in (-\infty,0)$ and the underdamped regime $\lambda \in (0,1)$ for any $\mu>0$, as well as the borderline critical case $\lambda=1$ under the improved condition $\mu>2\varepsilon^2$.

\end{abstract}
\vspace{2mm}

\vspace{2mm}


\section{Introduction}
\subsection{Presentation of the model and aims of the paper}

Nonautonomous coefficients naturally arise in evolution equations posed on expanding backgrounds, that is, media whose spatial scale increases with time, where the expansion induces a friction-like term, commonly referred to in cosmology as \emph{Hubble friction} \cite{Ryden}. Related effects can also be observed in physical wave systems whose material properties vary in time \cite{2,1}. From a mathematical viewpoint, time-dependent dissipation has been extensively studied for wave equations and related hyperbolic models.

Motivated by these considerations, we focus here on a prototypical hyperbolic system with nonautonomous damping coefficients, namely the compressible Euler equations with time-dependent damping in $\mathbb{R}^d$ ($d \geq 1$), which take the form
\begin{equation}\label{euler}
\left\{
\begin{aligned}
&\partial_{t}\rho^{\var}+\dive (\rho^{\var} u^{\var})=0,\\
&\var^{2} \partial_{t}(\rho^{\var}u^{\var})+\var^{2} \dive (\rho^{\var}u^{\var}\otimes u^{\var}) +\nabla P(\rho^{\var})+\frac{\mu\rho^{\var}u^{\var}}{(1+t)^{\lambda}}=0,
\end{aligned}
\right.
\end{equation}
where $\rho^\var=\rho^\var(t,x)\geq0$ is the density, $u^\var=u^\var(t,x)\in\mathbb{R}^{d}$ is the velocity, $\varepsilon$ is the relaxation parameter, the time-dependent friction coefficient takes the form $\frac{\mu}{(1+t)^{\lambda}}$ with $\lambda\leqslant1$ and $\mu>0$, and the pressure function $P(\rho)$ is assumed to satisfy 
\begin{align}
    &P(\rho)\in C^{\infty}(\mathbb{R}_{+})\quad \text{and}\quad P'(\rho)>0~~\text{for}~~\rho>0.\label{P}
\end{align}
We consider the Cauchy problem for \eqref{euler} supplemented with the initial data
\begin{align}
(\rho^\var, u^\var)(0,x)=(\rho^\var_0, u_0^\var)(x).\label{d}
\end{align}
Concerning the friction coefficient, we will distinguish the following two cases:
\begin{itemize}
\item Case $\lambda<0$: the damping effect is stronger than in the constant-coefficient case; this regime is referred to as the \emph{overdamped} case.

\item Case $\lambda>0$: the damping effect is weaker than in the constant-coefficient case; this regime is referred to as the \emph{underdamped} case. When $0<\lambda<1$, global stability is expected for small perturbations, whereas for $\lambda\geq 1$, solutions may blow up in finite time. In particular, $\lambda=1$ is called the {\emph{critical}} case, and whether blow-up occurs depends on the value of $\mu$.
\end{itemize}
For related classifications, we refer to \cite{ji1,ji2} and the references therein.

If $(\rho^\var,u^\var)$ is a global solution to \eqref{euler}--\eqref{d}, we formally denote
\begin{align*}
(\rho^*,u^*):=\lim_{\var\rightarrow0}(\rho^\var, u^\var),\qquad 
\rho_0^*=\lim_{\var\rightarrow 0} \rho_0^\var.
\end{align*}
As $\var\rightarrow 0$, one expects that the dynamics of system \eqref{euler} are governed by a porous-medium-type diffusion model with time-dependent coefficients:
\begin{equation}\label{pm}
\left\{
\begin{aligned}
 &\partial_t\rho^{*}-\frac{(1+t)^{\lambda}}{\mu}\Delta P(\rho^{*})=0,\\
 &\rho^*(0,x)=\rho_0^*(x),
\end{aligned}
\right.
\end{equation}
and that $u^*$ is determined by Darcy’s law
\begin{align}\label{darcy}
\rho^{*}u^{*}=-\frac{(1+t)^{\lambda}}{\mu}\nabla P(\rho^{*}).
\end{align}

In this paper, we address two main objectives: the global-in-time well-posedness of \eqref{euler} in critical Besov spaces for perturbations of constant equilibria, and the rigorous justification of the relaxation limit from \eqref{euler} to the diffusive system \eqref{pm}--\eqref{darcy} as $\varepsilon\to0$. Theorem \ref{theorem11} establishes the first objective and also provides estimates uniform in $\varepsilon$. Based on these uniform estimates, Theorem \ref{theorem13} shows the second objective, namely, the uniform-in-time strong convergence of the relaxation limit with an explicit convergence rate.

\subsection{Literature review}
\label{sec:review}

\subsubsection{Link with the nonlinear wave equation}

System \eqref{euler} can be rewritten as a nonlinear wave equation with time-dependent damping, which naturally connects our problem to the classical diffusion phenomenon for damped wave equations. Indeed, in the case $\var=1$, when a solution to \eqref{euler} is viewed as a perturbation of a non-vacuum constant equilibrium $(\bar{\rho},0)$, the wave formulation of the linearized equation for the modified density perturbation  (see \eqref{eulerre}) reads as follows:
\begin{equation}\label{linearwave}
\partial_t^2 n-P'(\bar{\rho})\Delta n+\frac{\mu}{(1+t)^\lambda}\partial_t n=0.
\end{equation}
According to the analysis of Wirth \cite{Wirth:2004MMAS,Wirth:2006JDE,Wirth:2007JDE}, when $\lambda<1$, the damping is effective and the large-time behavior is of parabolic type. In particular, solutions to \eqref{linearwave} are asymptotically related to those of a heat equation with time-dependent diffusion coefficients, namely, the linearized equation associated with \eqref{pm}. In contrast, when $\lambda>1$, the damping becomes ineffective, and the asymptotic behavior is governed instead by the free wave equation
\begin{equation*}
\partial_t^2 n-P'(\bar{\rho})\Delta n=0.
\end{equation*}

The case $\lambda=1$ is critical. In this regime, the decay properties of \eqref{linearwave} depend sensitively on the value of $\mu$; see Wirth \cite{Wirth:2004MMAS}. {{In particular, Wirth showed that $\mu=2$ is the critical value for the sharp $L^2$--$L^2$ estimates of the associated energy operator: below this threshold, the decay order depends on $\mu$, whereas above it, increasing $\mu$ no longer improves the decay order.}} Interestingly, the same threshold also appears in the nonlinear compressible Euler equations. In one space dimension, Pan \cite{pan2} proved global existence for sufficiently small smooth perturbations when $\mu>2$, whereas Pan \cite{pan1} established finite-time blow-up for a class of small, compactly supported initial data satisfying suitable sign conditions when $0\leq\mu\leq2$.

In the critical case $\lambda=1$, there is a discrepancy between the nonlinear compressible Euler system and its diffusive asymptotics: the known one-dimensional global small-data theory for the Euler system requires $\mu>2$, whereas the diffusive profile \eqref{pm} remains stable for any $\mu>0$. Our results show that, in the presence of the relaxation parameter $\var$, this threshold is relaxed to $\mu>2\var^2$, so that as $\var\to0$, the condition on $\mu$ becomes consistent with the limiting diffusive equation.

\subsubsection{Literature review for the case $\varepsilon = 1$}
When $\var = 1$, system \eqref{euler} has been the subject of extensive investigations in the literature. 
In the constant damping coefficient  case $(\lambda = 0)$, system \eqref{euler} reduces to the well-known compressible Euler system with damping. Hsiao and Liu~\cite{HsiaoL:1992CMP} first observed that solutions to the one–dimensional damped Euler equations asymptotically approach the self–similar profile of the corresponding nonlinear porous–medium equation, the so–called \emph{diffusion wave}. 
Subsequently, Nishihara~\cite{Nishihara:1996JDE}, Nishihara \emph{et al.}~\cite{NishiharaWY:2000JDE}, and Mei~\cite{Mei:2010SIAM} quantified the convergence rates toward such diffusion waves in various functional frameworks. 
Sideris \emph{et al.}~\cite{si1} established the global existence and time–decay of small–amplitude smooth solutions near a non–vacuum constant state in three dimensions. 
Tan and Wu~\cite{TanW:2012JDE} as well as Tan and Wang~\cite{TanW:2013JDE} further improved these results by employing the Besov–space approach, and the optimal pointwise decay for multidimensional systems was later derived by Wang and Yang~\cite{wang1}. 
For initial data containing vacuum, the existence of entropy solutions and their $L^1$–weak convergence toward the Barenblatt self–similar profile were obtained in a series of works by Huang and Pan~\cite{HuangP:2003ARMA}, Huang \emph{et al.}~\cite{HuangMP:2005ARMA,HuangPW:2011ARMA}, and Geng and Huang~\cite{geng0}. 
Moreover, the convergence of classical solutions in the \emph{physical-vacuum} regime was further analyzed by Luo and Zeng~\cite{LuoZ:2016CPAM} and by Zeng~\cite{Zeng:2017ARMA,Zeng:2021ARMA}. Then, extensions of global existence and time-decay rates near equilibrium to critical spaces were obtained in non-homogeneous settings by Kawashima and Xu in \cite{XK1,XK2} and in some hybrid homogeneous settings by Crin-Barat and Danchin in \cite{c1,c2,c3}. Recently, the optimality of decay rates and regularity criteria have been studied by Shou, Xu and Zhang in \cite{SXZ}.

When the damping coefficient depends on time (that is, $\lambda\neq 0$), the dynamics of the compressible Euler system become more delicate.
In the one–dimensional case, Pan~\cite{pan1,pan2} proved that if $\lambda \in (0,1)$ with $\mu > 0$ or $\lambda = 1$ with $\mu > 2$, and the initial data are small, smooth perturbations of a non–vacuum constant state, then the corresponding classical solution exists globally in time. 
Chen \emph{et al.}~\cite{chen1} subsequently extended the global existence result to certain classes of large initial data. 
When $\lambda > 1,\, \mu > 0$ or $\lambda = 1,\, \mu \le 2$, the $C^1$ solution blows up in finite time; the blow–up mechanism was investigated by Sugiyama~\cite{su1}. 
Convergence toward the diffusion wave profile was independently established by Cui \emph{et al.}~\cite{CuiYZZ} and Li \emph{et al.}~\cite{lht1,lht2}, where the asymptotic states at spatial infinity are distinct. 
In the critical damping case, Geng \emph{et al.}~\cite{geng1} further proved convergence toward the asymptotic profile with an explicit rate depending on the physical parameter $\mu$.

For higher-dimensional cases $d=2,3$, Hou and Yin~\cite{hou1} first established global smooth solutions when $\lambda\in(0,1)$ with $\mu>0$, or when $\lambda=1$ with $\mu>3-d$, provided that the initial perturbation is small, smooth, compactly supported, and curl-free. In the three-dimensional potential-flow setting, Hou, Witt, and Yin~\cite{hou2} established corresponding global-existence and blow-up results. In contrast, when $\lambda > 1,\, \mu > 0$ or $\lambda = 1,\, \mu \le 3 - d$, the solution blows up in finite time. 
The decay rates for multidimensional solutions in the range $\lambda \in (0,1)$ were first obtained by Pan~\cite{pan3} and later refined by Ji and Mei~\cite{ji1,ji2}. 
The $L^1$–weak convergence to the generalized Barenblatt self–similar solution was established by Geng \emph{et al.}~\cite{geng2}, while the strong convergence in the physical–vacuum regime toward the generalized Barenblatt profile was rigorously justified by Pan~\cite{pan4,pan5} in the one–dimensional and spherically symmetric three–dimensional settings.

The regularity index $d/2+1$ is critical since $\dot{\B}^{d/2+1}_{2,1}$ is continuously embedded in the space of globally Lipschitz functions, and the control of the Lipschitz norm is well known to be crucial for preventing finite-time blow-up in hyperbolic systems; see, for instance, \cite{Dafermos1}. We also refer to \cite{IPS} for the ill-posedness of hyperbolic systems in $H^s$ with $s<d/2+1$. To ensure both the Banach space structure and the $C^1$ regularity, a natural minimal regularity space for the initial data is  $\dot{\B}^{d/2}_{2,1}\cap\dot{\B}^{d/2+1}_{2,1} $. 

While the mathematical results in critical spaces are now rather well understood for the constant-damping case, the corresponding existence theory for compressible Euler equations with time-dependent damping remains much less complete.

\subsubsection{Case $\varepsilon>0$}

The study of relaxation limits for hyperbolic relaxation systems has a long history. 
The pioneering results in one spatial dimension are due to Marcati, Milani, and Secchi~\cite{MMS}, who employed the method of compensated compactness. 
Further contributions were made by Liu~\cite{Liu}, Marcati and Milani~\cite{MarcatiM:1990JDE}, as well as Marcati and Rubino~\cite{marcatiparabolicrelaxation}, who developed a complete hyperbolic–to–parabolic relaxation theory in one dimension.
For the isothermal Euler equations, Junca and Rascle~\cite{Junca} established convergence to the heat equation for large $BV$ data away from vacuum.  
In several dimensions, the uniform regularity estimates and  the weak relaxation limit for the damped Euler system to the porous medium equation have been proved in \cite{CoulombelGoudon,XuWang}. Concerning the explicit convergence rates, the first author and Danchin \cite{c3} developed a frequency-localized functional setting to derive the strong relaxation limit with explicit convergence rates for ill-prepared data. The functional techniques have been adapted to some singular limits for different models with non-standard dissipation structures (see \cite{chs:siam,cst,chi1,csJinXin}). By using direct error estimates in Sobolev spaces without frequency localization,  Crin-Barat, Peng and Shou \cite{CBPS} obtained global convergence rates for global solutions in the ill-prepared setting.

Despite these advances, the rigorous justification of such relaxation limits has so far been established only for \emph{constant damping coefficients} (that is, $\lambda=0$).

In this paper, we investigate the relaxation limit problem in the time-dependent damping case $\lambda\neq0$, in which dissipation has a non-autonomous structure. The main challenge is that, in order to recover the full dissipation, the analysis must involve time-dependent weights, whose time derivatives are, however, difficult to handle. One of our main contributions is to control these terms by carefully exploiting the interplay between time and frequency.

\subsection{Outline of the paper}

The rest of the paper is organized as follows. In Section \ref{sec:spect}, we present a formal spectral analysis involving the parameter $\var$ and introduce the functional framework used throughout the paper. Section \ref{sec:mainresult} is devoted to the statement of the main results, namely, the global well-posedness result (Theorem \ref{theorem11}) and the justification of the relaxation limit (Theorem \ref{theorem13}). The proof of Theorem \ref{theorem11} is given in Section \ref{sec:globalex}, while Section \ref{sec:relaxation} is devoted to the proof of Theorems \ref{theorem12} and \ref{theorem13}. Some technical lemmas are collected in the Appendix.

\section{Spectral analysis and functional framework}\label{sec:spect}
\subsection{Spectral analysis involving the relaxation parameter}

Under the condition \eqref{P}, if $\rho^\var$ is a small perturbation of $\bar{\rho}$, we introduce the unknowns
\begin{equation*}   n^\var:=\int_{\bar{\rho}}^{\rho^\var}\frac{P'(s)}{s}\,ds\quad\quad\text{and}\quad n_{0}^\var:=\int_{\bar{\rho}}^{\rho^\var_{0}}\frac{P'(s)}{s}\,ds.
\end{equation*}
The Cauchy problem of System \eqref{euler} with the initial data $(\rho_{0}^\var,u_{0}^\var)$ can be reformulated as
\begin{equation}\label{eulerre}
\left\{
    \begin{aligned}
    &\partial_{t}n^\var+u^\var\cdot \nabla n^\var+(P'(\bar{\rho})+G(n^\var))\dive u^\var=0,\\
    &\var^2(\partial_{t}u^\var+u^\var\cdot \nabla u^\var)+\nabla n^\var+\frac{\mu}{(1+ t)^{\lambda}} u^\var=0,\\
    &(n^\var,u^\var)(0,x)=(n_{0}^\var,u_{0}^\var)(x),
    \end{aligned}
    \right.
\end{equation}
with the nonlinear term
$$
G(n^\var):=P'(\rho^\var)-P'(\bar{\rho}).
$$
 Since $P$ is a smooth function, we observe that $G(n^\var)$ also depends smoothly on $n^\var$, and  $n^\var$ and $\rho^\var-\bar\rho$ are equivalent variables in a neighborhood of the equilibrium state. To simplify the notation, we shall drop the superscript $\var$ in the following analysis.

In order to understand the behavior of solutions to \eqref{eulerre} with respect to the time-dependent friction coefficient, we perform a spectral analysis of the linearized system. In the Hodge decomposition, we denote the compressible part $m= \var\Lambda^{-1}\dive u$ and the incompressible part $\omega=\var\Lambda^{-1}\nabla\times u$ with $\Lambda^{\sigma}:=\mathcal{F}^{-1}(|\xi|^{\sigma}\mathcal{F}(\cdot))$. In dimension one, the incompressible component is absent. Thus, the linearization of system \eqref{eulerre} reads:
 \begin{equation}\nonumber
\begin{aligned}
& \partial_{t}
\left(\begin{matrix}
   n  \\
   m\\
  \end{matrix}\right)
  =\mathbb{A}\left(\begin{matrix}
  n \\
  m
     \end{matrix}\right),\quad  \mathbb{A}:=\left(\begin{matrix}
0.                             &      - \frac{1}{\var}P'(\bar{\rho})\Lambda\\
  \frac{1}{\var}\Lambda                      &  -\frac{1}{\var^2 b( t)}
  \end{matrix}\right),
 \quad\quad \partial_{t} \omega +\frac{1}{ \var^2 b( t)}\omega=0,
   \end{aligned}
\end{equation}
where \begin{align*}
    b(t)=\frac{(1+t)^\lambda}{\mu}.
\end{align*}
The eigenvalues of $\widehat{\mathbb{A}}(\xi)$, which is the Fourier symbol of the operator $\mathbb{A}$ defined previously, satisfy
\begin{equation}\nonumber
\begin{aligned}
&\lambda_{\pm}=-\frac{1}{2 \var^2 b( t)}\pm \frac{1}{2\var}\sqrt{ \frac{1}{\var^2 b^2( t)}-4P'(\bar{\rho})|\xi|^2}.
\end{aligned}
\end{equation}
\begin{itemize}
    \item In the low-frequency regime $|\xi|\ll\frac{1}{ \var b(t)}$, all the eigenvalues are real, and we have $\lambda_{+}\sim -b(t){ P'(\bar{\rho})}|\xi|^2$ and $\lambda_{-}\sim -\frac{1}{\var^2 b(t)} $.
    
        \item In the high-frequency regime $|\xi|\gg \frac{1}{ \var b(t)} $, the eigenvalues $\lambda_{\pm}$ are conjugate complex numbers and satisfy $\lambda_{\pm}\sim -\frac{1}{2\var^2 b(t)}\pm \frac{\sqrt{P'(\bar{\rho})}}{\var}|\xi| {\rm i}.$
\end{itemize}

In order to capture the optimal dissipation structures in each frequency regime, the above spectral analysis suggests separating the frequency space into two parts.  Since we shall use the Littlewood-Paley decomposition, we will in fact use the dyadic threshold $2^{J_{t}}\sim  \frac{1}{\var b(t)}$. Precisely, for $t\geqslant 0$, we set the threshold
\begin{align}
    &J_{t}:=\Big[\log_{2}\frac{1}{\var b(t)}\Big]-k_{0},         \label{J}
\end{align}
for some generic constant $k_{0}\in \Z$ and where $[a]$ is the integer part of a real number $a$. We find that in the autonomous coefficient case, i.e. $\lambda=0$, we have 
\begin{align*}
\quad\quad J_t\equiv J_0= \left[\log_2 \frac{\mu}{\var}\right]-k_0.
\end{align*}

\begin{figure}[!h]
    \centering
    \begin{minipage}[t]{0.48\textwidth}
        \centering
        \begin{tikzpicture}[xscale=0.85,yscale=0.85, thick]
            % axes
            \draw[-latex] (-0.2,0) -- (6.9,0);
            \draw[-latex] (0,-3.4) -- (0,3.4);

            % axis labels
            \node[left] at (0,3.2) {$j$};
            \node[below left] at (0,0) {$0$};
            \node[below] at (6.55,0) {$\infty$};

            % maximal time t
            \def\tmax{5.0}
            \pgfmathsetmacro{\Jt}{2.4-1.9*ln(1+\tmax)}

            % low-frequency shaded region, only up to t
            \fill[gray!12,opacity=0.32]
                plot[smooth,domain=0.03:\tmax,samples=120] (\x,{2.4-1.9*ln(1+\x)})
                -- (\tmax,-3.1) -- (0.03,-3.1) -- cycle;

            % threshold curve: ends at (t,J_t)
            \draw[line width=1.15pt]
                plot[smooth,domain=0:\tmax,samples=120] (\x,{2.4-1.9*ln(1+\x)});

            % dashed guides
            \draw[dashed] (\tmax,-3.1) -- (\tmax,3.1);
            \draw[dashed] (0,\Jt) -- (\tmax,\Jt);
            \fill (\tmax,\Jt) circle (1.5pt);

            % labels
            \node[below right] at (\tmax,0) {$t$};
            \node[left] at (0,\Jt) {$J_t$};

            % frequency-region labels
            \node at (2.3,-1.35) {Low frequencies};
            \node at (2.1,2.1) {High frequencies};
        \end{tikzpicture}

        \vspace{0.3em}
        \small (a) Underdamped case $0<\lambda\leq 1$
    \end{minipage}
    \hfill
    \begin{minipage}[t]{0.48\textwidth}
        \centering
        \begin{tikzpicture}[xscale=0.85,yscale=0.85, thick]
            % axes
            \draw[-latex] (-0.2,0) -- (6.9,0);
            \draw[-latex] (0,-3.4) -- (0,3.4);

            % axis labels
            \node[left] at (0,3.2) {$j$};
            \node[below left] at (0,0) {$0$};
            \node[below] at (6.55,0) {$\infty$};

            % maximal time t
            \def\tmax{5.0}
            \pgfmathsetmacro{\Jt}{-2.2+1.45*ln(1+\tmax)}

            % low-frequency shaded region, only up to t
            \fill[gray!12,opacity=0.32]
                plot[smooth,domain=0.03:\tmax,samples=120] (\x,{-2.2+1.45*ln(1+\x)})
                -- (\tmax,-3.1) -- (0.03,-3.1) -- cycle;

            % threshold curve: ends at (t,J_t)
            \draw[line width=1.15pt]
                plot[smooth,domain=0:\tmax,samples=120] (\x,{-2.2+1.45*ln(1+\x)});

            % dashed guides
            \draw[dashed] (\tmax,-3.1) -- (\tmax,3.1);
            \draw[dashed] (0,\Jt) -- (\tmax,\Jt);
            \fill (\tmax,\Jt) circle (1.5pt);

            % labels
            \node[below right] at (\tmax,0) {$t$};
            \node[left] at (0,\Jt) {$J_t$};

            % frequency-region labels
            \node at (2.2,-2.1) {Low frequencies};
            \node at (2.2,1.8) {High frequencies};
        \end{tikzpicture}

        \vspace{0.3em}
        \small (b) Overdamped case $\lambda<0$
    \end{minipage}
    \caption{Frequency-time partition}
    \label{fig:frequency-splitting}
\end{figure}
 In the case $\lambda\ne0$, the frequency threshold $J_t$ will also depend on {\emph{the time}}; this introduces substantial technical difficulties compared with the constant-damping case. Indeed, for a fixed frequency indexed by $j$, whether the mode belongs to the low- or high-frequency regime depends on time. To quantify this interplay between time and frequency, we define a time threshold $t_j$ as follows:
\begin{equation*}
t_j:=\max\left\{\left(\frac{\mu}{\var 2^{k_0+j}}\right)^{\frac{1}{\lambda}}-1,0\right\}, 
\end{equation*}
that is
\begin{equation*}
\text{if }\quad  \lambda\in(0,1],\quad t_j=\left\{
\begin{aligned}
& \left(\frac{\mu}{\var 2^{k_0+j}}\right)^{\frac{1}{\lambda}}-1,\quad  j\leq J_0,\\
& 0, \quad \quad \quad \quad\quad \,\,  \quad \quad  j\geq J_0+1,
\end{aligned}
\right. 
\end{equation*}
and
\begin{equation*}
\text{if} \quad \lambda<0,\quad  t_j=\left\{
\begin{aligned}
& \left(\frac{\mu}{\var 2^{k_0+j}}\right)^{\frac{1}{\lambda}}-1, \quad  j\geq J_0+1,\\
& 0, \quad \quad \quad \quad\,\, \quad \quad \quad  j\leq J_0.
\end{aligned}
\right. 
\end{equation*}
Roughly speaking, $t_j$ is the time such that $j=J_{t_j}$. More precisely, when $t_j>0$, we have 
\begin{equation*}2^{j}= (\var b(t_j))^{-1}2^{-k_{0}}, \end{equation*}
which means exactly that $j=J_{t_j}$. And if $t_j=0$, it means that in the underdamped case, $0<\lambda\leqslant 1$, for all $j\geq J_0+1$, the corresponding modes always remain in the high-frequency regime, whereas in the overdamped case, the opposite occurs.

\subsection{Functional framework}
We define the Besov semi-norms adapted to the spectral analysis described in the previous section. For a fixed time $\tau>0$, we define 
\begin{align*}
    \|u(\tau)\|_{\dot{\B}_{p,r}^s}^{\ell}:=\|\{2^{js}\|\dot{\Delta}_j u\|_{L^p}\}_{j\leqslant J_\tau}\|_{\ell^r}\quad{\rm and} \quad\|u(\tau)\|_{\dot{\B}_{p,r}^s}^{h}:=\|\{2^{js}\|\dot{\Delta}_j u\|_{L^p}\}_{j\geqslant J_\tau+1}\|_{\ell^r},
\end{align*}
where the dyadic Littlewood-Paley blocks are the classical ones defined in \cite[Section 2]{bahouri1} (see also the appendix).
For fixed $j$ and $t$, we define
\begin{align*}
    I_{j,t}^\ell:=\{0\leq\tau\leq t\mid j\leqslant J_\tau\}\quad{\rm and} \quad I_{j,t}^h:=\{0\leq\tau\leq t\mid j\geqslant J_\tau+1\}.
\end{align*}

\begin{remark}\label{remark11}
Since $b(t)$ is monotone, both $I_{j,t}^\ell$ and $I_{j,t}^h$ are intervals. More precisely, according to the position of $t_j$ relative to $[0,t]$, one has the following descriptions:
\begin{itemize}
\item If $0<\lambda\leq 1$, then
\[
I_{j,t}^\ell=
\begin{cases}
[0,t], & t_j\ge t,\\
[0,t_j], & t_j\in[0,t],
\end{cases}
\qquad
I_{j,t}^h=
\begin{cases}
\emptyset, & t_j\ge t,\\
[t_j,t], & t_j\in[0,t],
\end{cases}
\]

\item If $\lambda<0$, then
\[
I_{j,t}^\ell=
\begin{cases}
\emptyset, & t_j\ge t,\\
[t_j,t], & t_j\in[0,t],
\end{cases}
\qquad
I_{j,t}^h=
\begin{cases}
[0,t], & t_j\ge t,\\
[0,t_j], & t_j\in[0,t].
\end{cases}
\]
\end{itemize}
\end{remark}

\begin{remark}
From Figures \ref{fig:frequency-splitting} and \ref{fig:time-splitting}, we see that $I_{j,t}^\ell$ and $I_{j,t}^h$ simply denote the time sets during which the dyadic block $\ddj$ localizes in the low- and high-frequency regimes, respectively.
\end{remark}

\begin{figure}[!h]
    \centering

    \begin{minipage}{0.95\textwidth}
        \centering
        \begin{tikzpicture}[xscale=0.48,yscale=0.48, thick]
            % axis
            \draw[-latex] (-1,0) -- (21,0);
            \draw (21,0) node[below] {$t$};

            % ticks
            \draw (-1,0) node {$|$};
            \draw (10,0) node {$|$};

            \draw (-1,-0.25) node[below] {$0$};
            \draw (10,-0.25) node[below] {$t_j$};

            % labels
            \draw (4.5,0.0) node[above] {Low frequencies};
            \draw (4.5,-0.75) node {$\dot{\B}^{\frac{d}{2}}_{2,1}$};

            \draw (15.5,0.0) node[above] {High frequencies};
            \draw (15.5,-0.75) node {$\dot{\B}^{\frac{d}{2}+1}_{2,1}$};
        \end{tikzpicture}

        \vspace{0.3em}
        \small (a) Underdamped case $0<\lambda\leq 1$
    \end{minipage}

    \vspace{0.8em}

    \begin{minipage}{0.95\textwidth}
        \centering
        \begin{tikzpicture}[xscale=0.48,yscale=0.48, thick]
            % axis
            \draw[-latex] (-1,0) -- (21,0);
            \draw (21,0) node[below] {$t$};

            % ticks
            \draw (-1,0) node {$|$};
            \draw (10,0) node {$|$};

            \draw (-1,-0.25) node[below] {$0$};
            \draw (10,-0.25) node[below] {$t_j$};

            % labels
            \draw (4.5,0.0) node[above] {High frequencies};
            \draw (4.5,-0.75) node {$\dot{\B}^{\frac{d}{2}+1}_{2,1}$};

            \draw (15.5,0.0) node[above] {Low frequencies};
            \draw (15.5,-0.75) node {$\dot{\B}^{\frac{d}{2}}_{2,1}$};
        \end{tikzpicture}

        \vspace{0.3em}
        \small (b) Overdamped case $\lambda<0$
    \end{minipage}

    \caption{Frequency-dependent time splitting}
    \label{fig:time-splitting}
\end{figure}

For $\varrho\geqslant1$, we adapt the Chemin-Lerner spaces to more general spaces involving time-dependent coefficients:
\begin{align*}
\text{For }0<\lambda\leq1,\qquad
\|u\|_{\widetilde{L}^{\varrho}_{t}(\dot{\B}^{s}_{p,1})}^{\ell}
&:=
\sum_{\substack{j\in\Z\\ t_j>0}}
2^{js}
\left(
\int_{I_{j,t}^{\ell}}
\|\ddj u(\tau)\|_{L^p}^{\varrho}\,d\tau
\right)^{\frac1\varrho},
\\
\phantom{\text{For }0<\lambda\leq1,\qquad}
\|u\|_{\widetilde{L}^{\varrho}_{t}(\dot{\B}^{s}_{p,1})}^{h}
&:=
\sum_{\substack{j\in\Z\\ t_j<t}}
2^{js}
\left(
\int_{I_{j,t}^{h}}
\|\ddj u(\tau)\|_{L^p}^{\varrho}\,d\tau
\right)^{\frac1\varrho},
\\[0.3em]
\text{For }\lambda<0,\qquad
\|u\|_{\widetilde{L}^{\varrho}_{t}(\dot{\B}^{s}_{p,1})}^{\ell}
&:=
\sum_{\substack{j\in\Z\\ t_j<t}}
2^{js}
\left(
\int_{I_{j,t}^{\ell}}
\|\ddj u(\tau)\|_{L^p}^{\varrho}\,d\tau
\right)^{\frac1\varrho},
\\
\phantom{\text{For }\lambda<0,\qquad}
\|u\|_{\widetilde{L}^{\varrho}_{t}(\dot{\B}^{s}_{p,1})}^{h}
&:=
\sum_{\substack{j\in\Z\\ t_j>0}}
2^{js}
\left(
\int_{I_{j,t}^{h}}
\|\ddj u(\tau)\|_{L^p}^{\varrho}\,d\tau
\right)^{\frac1\varrho},
\\[0.3em]
\text{For }\lambda=0,\qquad
\|u\|_{\widetilde{L}^{\varrho}_{t}(\dot{\B}^{s}_{p,1})}^{\ell}
&:=
\sum_{\substack{j\in\Z\\ j\leq J_0}}
2^{js}
\left(
\int_0^t
\|\ddj u(\tau)\|_{L^p}^{\varrho}\,d\tau
\right)^{\frac1\varrho},
\\
\phantom{\text{For }\lambda=0,\qquad}
\|u\|_{\widetilde{L}^{\varrho}_{t}(\dot{\B}^{s}_{p,1})}^{h}
&:=
\sum_{\substack{j\in\Z\\ j\geq J_0+1}}
2^{js}
\left(
\int_0^t
\|\ddj u(\tau)\|_{L^p}^{\varrho}\,d\tau
\right)^{\frac1\varrho}.
\end{align*}
where, for $\varrho=+\infty$, the usual convention (involving the essential supremum $\sup_{[a,b]}f(\tau)$) is adopted.  The summation ranges in the above definitions depend on the monotonicity of the frequency threshold $J_t$. If $0<\lambda\leq1$, then $J_t$ decreases with time, so each dyadic block may pass from the low-frequency region to the high-frequency region at $t_j$. Hence, the low-frequency norm is summed over $t_j>0$, while the high-frequency norm is summed over $t_j<t$. If $\lambda<0$, then $J_t$ increases with time and the transition is reversed, so the two summation conditions must be interchanged. With these definitions, each dyadic block is integrated exactly over the time interval on which it belongs to the corresponding frequency region. By Fubini's Theorem, we observe that, for $\varrho=1$,
\begin{align*}
 \text{For } 0<\lambda\leq 1,&\quad   \|u\|_{\widetilde{L}^{1}_{t}(\dot{\B}^{s}_{p,1})}^{\ell}=\sum_{\substack{j\in\Z\\ t_j>0}}2^{js}\int_{I_{j,t}^\ell}\|\ddj u(\tau)\|_{L^p}\,d\tau=\int_0^t\|u(\tau)\|_{\dot{\B}^{s}_{p,1}}^{\ell}\,d\tau,\\
 & \quad   \|u\|_{\widetilde{L}^{1}_{t}(\dot{\B}^{s}_{p,1})}^{h}=\sum_{\substack{j\in\Z\\t_j<t}}2^{js}\int_{I_{j,t}^h}\|\ddj u(\tau)\|_{L^p}\,d\tau=\int_{0}^t\|u(\tau)\|_{\dot{\B}^{s}_{p,1}}^{h}\,d\tau,
 \\
  \text{for } \lambda<0,&\quad   \|u\|_{\widetilde{L}^{1}_{t}(\dot{\B}^{s}_{p,1})}^{\ell}=\sum_{\substack{j\in\Z\\ t_j<t}}2^{js}\int_{I_{j,t}^\ell}\|\ddj u(\tau)\|_{L^p}\,d\tau=\int_{0}^t\|u(\tau)\|_{\dot{\B}^{s}_{p,1}}^{\ell}\,d\tau,\\
  &\quad   \|u\|_{\widetilde{L}^{1}_{t}(\dot{\B}^{s}_{p,1})}^{h}=\sum_{\substack{j\in\Z\\ t_j>0}}2^{js}\int_{I_{j,t}^h}\|\ddj u(\tau)\|_{L^p}\,d\tau=\int_0^t\|u(\tau)\|_{\dot{\B}^{s}_{p,1}}^{h}\,d\tau.
\end{align*}
This shows that the $L^1$-in-time Chemin–Lerner spaces coincide with the usual Lebesgue–Besov spaces. On the other hand, by Minkowski's inequality, the $L^\infty$-in-time Chemin–Lerner norms are stronger than the corresponding $L^\infty$-in-time Lebesgue–Besov norms.

Classical approaches to \eqref{eulerre} consist of reformulating the system as a second–order wave equation with time-dependent coefficients and then
applying the analytic tools available for wave equations (cf. \cite{pan2,hou1}). In contrast, in the present work, we develop a direct hypocoercive energy method for the first–order
hyperbolic system \eqref{eulerre}, without passing through the wave formulation.

A key difficulty of the problem is that, owing to the time-dependent coefficient $b(t)$, the frequency threshold is no longer fixed, so a given mode may exhibit different spectral behaviors at different times. As a consequence, the standard energy method based on a fixed low/high-frequency decomposition, which is well adapted to the constant-coefficient setting (see, e.g., \cite{c3}), is no longer directly applicable. In addition, compared with the constant-coefficient case, an additional difficulty is caused by the extra $b'(t)$-terms arising in the energy-dissipation estimates.

Our analysis is based on a refined frequency decomposition and a careful low/high-frequency analysis within the Littlewood–Paley framework. More precisely, the low-frequency part is estimated in the critical space $\dot{\B}^{d/2}_{2,1}$, while the high-frequency part is treated in $\dot{\B}^{d/2+1}_{2,1}$. For each dyadic block, the energy estimates must be performed separately on the time intervals $I_{j,t}^\ell$ and $I_{j,t}^h$. Moreover, once the analysis is split at the transition time $t_j$, the estimates on these two intervals have to be patched together rather than treated independently. As a consequence, the low- and high-frequency analyses are no longer decoupled (see Figures \ref{fig:frequency-splitting} and \ref{fig:time-splitting}). Furthermore, the additional $b'(t)$-terms arising in the energy estimates can be controlled by converting the $L^1$-in-time control into an $L^\infty$-in-time bound (see \eqref{b'1}) and using the structural conditions on $\lambda$ and $\mu$ (see \eqref{condHF}). This strategy allows us to exploit the maximal $L^{1}$-in-time integrability of the dissipation in a low-regularity critical Besov framework.

\section{Main results}\label{sec:mainresult}

Our first result concerns the global well-posedness of the Cauchy problem for System \eqref{euler} in the critical regularity setting and establishes uniform regularity estimates with respect to $\var$.

\begin{theorem}\label{theorem11}  Let $d\geq1$, $-\infty<\lambda\leq 1$, $\varepsilon\in(0,1]$, $\bar\rho>0$ and
\begin{equation*}
\begin{cases}
\mu >0, & \text{if } \,\, \lambda < 1, \\[0.3em]
\mu >2\varepsilon^2, & \text{if }\,\, \lambda = 1.
\end{cases}
\end{equation*}
There exists a constant $\delta_{0}>0,$ independent of $\varepsilon$, such that if the initial data $(\rho^\varepsilon_{0},u^\varepsilon_{0})$ satisfies $(\rho^\varepsilon_{0}-\bar{\rho},u^\varepsilon_{0})\in \dot{\B}^{\frac{d}{2}}_{2,1}\cap \dot{\B}^{\frac{d}{2}+1}_{2,1}$ and
\begin{equation*}
\|(\rho^\varepsilon_{0}-\bar{\rho},{{\varepsilon}}u^\varepsilon_{0})\|_{\dot{\B}^{\frac{d}{2}}_{2,1}}
+{{\varepsilon}}\|(\rho^\varepsilon_{0}-\bar{\rho},{{\varepsilon}}u^\varepsilon_{0})\|_{\dot{\B}^{\frac{d}{2}+1}_{2,1}}
\leq
\begin{cases}
\delta_0, & \text{if } \lambda<1,\\
\delta_0{{\theta_\varepsilon^4}}, & \text{if } \lambda=1.
\end{cases}
\end{equation*}
then the Cauchy problem \eqref{euler}-\eqref{d} admits a unique global classical solution $(\rho^\varepsilon,u^\varepsilon)$ that satisfies
\begin{equation*}
\begin{aligned}
(\rho^\varepsilon-\bar{\rho}, u^\varepsilon)\in \mathcal{C}(\mathbb{R}^+;\dot{\B}^{\frac{d}{2}}_{2,1}\cap \dot{\B}^{\frac{d}{2}+1}_{2,1})
\end{aligned}
\end{equation*}
and
\begin{equation*}
\begin{aligned}
&\|(\rho^\varepsilon-\bar{\rho},\varepsilon u^\varepsilon)\|_{\widetilde{L}^{\infty}_{t}(\dot{\B}^{\frac{d}{2}}_{2,1})}+\varepsilon \|(1+\tau)^{\lambda}(\rho^\varepsilon-\bar{\rho},\varepsilon u^\varepsilon)\|_{\widetilde{L}^{\infty}_{t}(\dot{\B}^{\frac{d}{2}+1}_{2,1})}\\
&\quad+\|(1+\tau)^{\lambda}(\rho^\varepsilon-\bar{\rho})\|_{L^{1}_{t}(\dot{\B}^{\frac{d}{2}+2}_{2,1})}^{\ell}+\frac{\theta_\varepsilon}{\varepsilon}\|\rho^\varepsilon-\bar{\rho}\|_{L^{1}_{t}(\dot{\B}^{\frac{d}{2}+1}_{2,1})}^{h}+\theta_\varepsilon\|u^\varepsilon\|_{L^1_{t}(\dot{\B}^{\frac{d}{2}+1}_{2,1})}\\
&\quad +\frac{\theta_\varepsilon}{\varepsilon}\Big\|\nabla P(\rho^\varepsilon)+\frac{\mu}{(1+\tau)^{\lambda}}\rho^\varepsilon u^\varepsilon \Big\|_{L^1_{t}(\dot{\B}^{\frac{d}{2}}_{2,1})}\\
\leq &C{{\theta_\varepsilon^{-\frac32}}}\Big(\|(\rho^\varepsilon_{0}-\bar{\rho},{{\varepsilon}} u^\varepsilon_{0})\|_{\dot{\B}^{\frac{d}{2}}_{2,1}}+{{\varepsilon}}\|(\rho^\varepsilon_{0}-\bar{\rho},{{\varepsilon}} u^\varepsilon_{0})\|_{\dot{\B}^{\frac{d}{2}+1}_{2,1}} \Big)\quad\text{for all}\quad t>0,
\end{aligned}
\end{equation*}
where $C>0$ is a generic constant, and $\theta_\varepsilon=1$ for $\lambda<1$, while $\theta_\varepsilon=1-\frac{2\varepsilon^2}{\mu}$ for $\lambda=1$.
\end{theorem}

\begin{remark}\normalfont
Theorem \ref{theorem11} provides the first result on global well-posedness of solutions to the compressible Euler equations with time-dependent damping in the critical regularity setting. It covers the overdamped case $\lambda<0$ and the underdamped case $0<\lambda\leq 1$. In contrast to earlier works, our approach is purely energy-based and avoids techniques tailored to time-dependent wave equations. 
\end{remark}

\begin{remark}\normalfont
In the critical borderline regime $\lambda=1$, Pan \cite{pan1,pan2} proved that, for the one-dimensional compressible Euler system with time-dependent damping (the relaxation parameter fixed to $\var=1$), classical solutions arising from small perturbations exist globally-in-time when $\mu>2$, whereas finite-time blow-up may occur when $0\leq \mu\leq 2$. 
In higher dimensions, the available global existence results in the critical case prior to this work are due to Hou and Yin~\cite{hou1} and Hou, Witt, and Yin~\cite{hou2}. These results require the additional irrotational constraint $\mathrm{curl}\,u_0=0$ and the condition $\mu>3-d$.

In the case $\varepsilon=1$, Theorem~\ref{theorem11} gives global existence and uniqueness in the critical case $\lambda=1$ \emph{without} imposing an irrotationality condition on the initial data. In general, we obtain the stability condition $\mu>2\var^2$ for every $d\geq1$. In particular, our analysis further yields global existence for \emph{all} $\mu>0$, provided $\var$ is sufficiently small.
\end{remark}

\begin{remark}\normalfont
In the critical case $\lambda=1$, the smallness condition involving $\theta_\var$ describes the degeneration of the high-frequency estimates as $\mu$ approaches 
$2\var^2$. This has no  effect on the relaxation limit, since for
any fixed $\mu>0$, as $\var\to0$, we have  $\mu-2\var^2\to\mu>0$ and $\theta_\var\to1$. 
\end{remark}

We next establish a global existence result for the porous medium equation \eqref{pm}. 
\begin{theorem}\label{theorem12}
Let $d\geq1$, $-\infty<\lambda\leq 1$, $\mu>0$ and $\bar{\rho}>0$. Let the initial data $\rho_0^*$ satisfy $\rho_0^*-\bar{\rho}\in \dot{\B}^{\frac{d}{2}}_{2,1}$. There exists a constant $\delta_0^*>0$ such that if
\begin{align*}
\|\rho_0^*-\bar{\rho}\|_{\dot{\B}^{\frac{d}{2}}_{2,1}}\leq \delta_0^*,
\end{align*}
then a unique global solution $\rho^*$ to the Cauchy problem \eqref{pm} exists, satisfies $\rho^*-\bar{\rho}\in \mathcal{C}(\mathbb{R}^+;\dot{\B}^{\frac{d}{2}}_{2,1})$ and
\begin{equation*}
\begin{aligned}
&\quad \|\rho^*-\bar{\rho}\|_{\widetilde{L}^{\infty}_t(\dot{\B}^{\frac{d}{2}}_{2,1})}+\|(1+\tau)^{\lambda}(\rho^*-\bar{\rho})\|_{L^1_t(\dot{\B}^{\frac{d}{2}+2}_{2,1})}+\|u^*\|_{L^1_t(\dot{\B}^{\frac{d}{2}+1}_{2,1})}\leq C\|\rho_0^*-\bar{\rho}\|_{\dot{\B}^{\frac{d}{2}}_{2,1}},
\end{aligned}
\end{equation*}
for all $t\geq0$. Here, $C>0$ is a generic constant, and $u^*$ is given by Darcy's law \eqref{darcy}.
\end{theorem}

Moreover, we justify the validity of the relaxation limit convergence and establish global-in-time error estimates between \eqref{euler}-\eqref{d} and \eqref{pm}-\eqref{darcy}. 
\begin{theorem}\label{theorem13}
{{Let $d\geq1$, $-\infty<\lambda\leq1$ with $\lambda\neq0$, $\bar\rho>0$, and $\mu>0$ be fixed. Assume that $\varepsilon$ satisfies $0<\varepsilon\leq1$ if $\lambda<1$ and is sufficiently small if $\lambda=1$.}} Let $(\rho^\var, u^\var)$ and $\rho^*$ be the global solutions to the problems \eqref{euler}-\eqref{d} and \eqref{pm} obtained in Theorems \ref{theorem11} and \ref{theorem12}, respectively, and let $u^*$ be given by Darcy's law \eqref{darcy}. 
\begin{itemize}

\item {\rm(Overdamped case)}: for $\lambda<0$,  there exists a uniform constant $C>0$ such that 
\begin{equation}\label{error1}
\begin{aligned}
&\|\rho^\var-\rho^*\|_{\widetilde{L}^{\infty}_t(\dot{\B}^{\frac{d}{2}-1}_{2,1})}+\|(1+\tau)^{\lambda}(\rho^\var-\rho^*)\|_{L^1_t(\dot{\B}^{\frac{d}{2}+1}_{2,1})}+\|u^\var-u^*\|_{L^1_t(\dot{\B}^{\frac{d}{2}}_{2,1})}\\
 \leq & C\|\rho_0^\var-\rho_0^*\|_{\dot{\B}^{\frac{d}{2}-1}_{2,1}}+C\var.
\end{aligned}
\end{equation}

\item {\rm(Underdamped case)}: for $0<\lambda\leq 1$, it holds for some uniform constant $C>0$ that
\begin{align}
&\|(1+\tau)^{-\lambda}(\rho^\var-\rho^*)\|_{\widetilde{L}^{\infty}_t(\dot{\B}^{\frac{d}{2}-1}_{2,1})}+\|\rho^\var-\rho^*\|_{L^1_t(\dot{\B}^{\frac{d}{2}+1}_{2,1})}+\|(1+\tau)^{-\lambda}(u^\var-u^*)\|_{L^1_t(\dot{\B}^{\frac{d}{2}}_{2,1})}\nonumber\\
\leq& C\|\rho_0^\var-\rho_0^*\|_{\dot{\B}^{\frac{d}{2}-1}_{2,1}}+C\var.\label{error2}
\end{align}
\end{itemize}
If $\|\rho_0^\var-\rho^*_0\|_{\dot{\B}^{\frac{d}{2}-1}_{2,1}}\leq \var^{q}$ {\rm(}$q>0${\rm)}, then the right-hand sides of \eqref{error1}-\eqref{error2} can be bounded by $\mathcal{O}(\var^{\min\{1,q\}})$. Consequently, as $\var\rightarrow0$, the solution of System \eqref{euler} converges strongly {\rm(}in the sense of \eqref{error1}-\eqref{error2}{\rm)} to the solution of \eqref{pm}-\eqref{darcy}.
\end{theorem}

\begin{remark}\normalfont
To the best of our knowledge, Theorem \ref{theorem13} provides the first rigorous justification of the singular limit from the compressible Euler system with time-dependent damping  toward the time-dependent porous medium equation and Darcy's law. 
\end{remark}

\begin{remark}\normalfont
We say that the initial data are well prepared if, as $\var\rightarrow0$, $u_0^\var=\mathcal{O}(1)$ holds and admits a limit  $\lim\limits_{\var\rightarrow0} u^\var_0$, and the compatibility condition (i.e., the convergence of $\eqref{euler}_2$ at $t=0$) holds, namely
\begin{align*}\nabla P(\rho^{\var}_0)+\mu\rho^{\var}_0 u^{\var}_0\rightarrow0\quad\text{as}\quad \var\rightarrow0.
\end{align*}
Otherwise, the data are said to be ill prepared.

Our analysis covers ill-prepared data and, in particular, allows singular initial velocities of size $u_0^\var=\mathcal{O}(\var^{-1})$.
\end{remark}

\begin{remark}\normalfont
Our analysis is done in $L^{2}$-based critical spaces. It is possible to extend our results to a $L^{p}$ framework in low frequencies when $\lambda \leqslant 0$ (see the constant damping case \cite{c3}). 
However, for the underdamped regime $0<\lambda \leq 1$, the  $L^{p}$ theory does not seem reachable with our current techniques. Indeed, 
the spectral analysis (and the classical work of Brenner \cite{Brenner1}) shows that we cannot expect $L^{p}$ estimates for the high-frequency part, essentially due to the presence of nontrivial imaginary parts in the eigenvalues. 
Moreover, since the frequency threshold separating low and high frequencies is time-dependent, every fixed frequency eventually enters the high-frequency region as time grows. Therefore, in the underdamped regime, it seems that one is restricted to $L^{2}$-based estimates at every frequency.
\end{remark}

\begin{remark}\normalfont
In a forthcoming study \cite{TSZ}, we aim to extend the present analysis to general partially dissipative systems satisfying the Shizuta-Kawashima condition. We also plan to treat more general time-dependent coefficients $b(t)$ and to identify the sharp conditions ensuring global-in-time existence and the large-time stability of the solutions. To this end, we build hypocoercive Lyapunov functionals in the spirit of Villani \cite{Villani}, further adapted to the hyperbolic setting by Beauchard and Zuazua \cite{BZ}.
\end{remark}

\section{Uniform global existence}\label{sec:globalex}

Throughout this section, we simplify the notations by omitting the superscript $\var$. To prove Theorem \ref{theorem11}, we first establish uniform {\emph{a priori}} estimates.  Define the energy functional
\begin{equation}\label{Xp}
\begin{aligned}
\mathcal{X}(t):&=\|(n,\var u)\|_{\widetilde{L}^{\infty}_{t}(\dot{\B}^{\frac{d}{2}}_{2,1})}^{\ell}+\|b(\tau)n\|_{L^{1}_{t}(\dot{\B}^{\frac{d}{2}+2}_{2,1})}^{\ell}+ \|u\|_{L^{1}_{t}(\dot{\B}^{\frac{d}{2}+1}_{2,1})}^{\ell}\\
&\quad+\var \|b(\tau)(n,\var u)\|_{\widetilde{L}^{\infty}_{t}(\dot{\B}^{\frac{d}{2}+1}_{2,1})}^{h}+\frac{\theta_\var}{\varepsilon}\|(n,\var u)\|_{L^{1}_{t}(\dot{\B}^{\frac{d}{2}+1}_{2,1})}^{h}\\
&\quad +\theta_\var^{\frac{1}{2}}\|b(\tau)^{-\frac{1}{2}}u\|_{\widetilde{L}^2_{t}(\dot{\B}^{\frac{d}{2}}_{2,1})}+\frac{\theta_\var}{\var}\|b(\tau)^{-1}u+\nabla n\|_{L^1_{t}(\dot{\B}^{\frac{d}{2}}_{2,1})},
\end{aligned}
\end{equation}
and the initial energy functional
\begin{equation}\label{Xp0}
\begin{aligned}
\mathcal{X}_{0}:=\|(n_{0},\var u_{0})\|_{\dot{\B}^{\frac{d}{2}}_{2,1}}^{\ell}+{\var}\|(n_{0},{\var}u_{0})\|_{\dot{\B}^{\frac{d}{2}+1}_{2,1}}^{h}.
\end{aligned}
\end{equation}
First, we observe that we have the following $L^{\infty}_t(L^\infty)$-control:
\begin{equation}
\begin{aligned}
\|(n,\var u)\|_{L^{\infty}_t(L^{\infty})}&\lesssim \|(n,\var u)\|_{L^{\infty}_t(\dot{\B}^{\frac{d}{2}}_{2,1})}\\
&\lesssim \|(n,\var u)\|_{\widetilde{L}^{\infty}_{t}(\dot{\B}^{\frac{d}{2}}_{2,1})}^{\ell}+\var 2^{k_0} \|b(\tau)(n,\var u)\|_{\widetilde{L}^{\infty}_{t}(\dot{\B}^{\frac{d}{2}+1}_{2,1})}^{h}\\
&\lesssim \mathcal{X}(t)\label{Linfty'}
\end{aligned}
\end{equation}
More generally, the following lemma will be used in our treatment of low–high frequency interactions.
\begin{lemma} \label{hltransfer}
For $f\in \mathcal{S}'(\R^d)$, $s\geqslant 0$ and $r\in [1,+\infty]$, we have
\begin{align}
     &\|f\|_{\widetilde{L}^{r}_{t}(\dot{\B}^{s}_{2,1})}
    \leqslant  \|f\|_{\widetilde{L}^{r}_{t}(\dot{\B}^{s}_{2,1})}^{\ell}+ 2^{k_0}\|\var b(\tau)f\|_{\widetilde{L}^{r}_{t}(\dot{\B}^{s+1}_{2,1})}^{h}, \label{lowhigh1}\\
    &\|f\|_{\widetilde{L}^{r}_{t}(\dot{\B}^{s}_{2,1})}\leqslant 2^{-k_0}\|(\var b(\tau))^{-1}f\|_{\widetilde{L}^{r}_{t}(\dot{\B}^{s-1}_{2,1})}^{\ell}+ \|f\|_{\widetilde{L}^{r}_{t}(\dot{\B}^{s}_{2,1})}^{h}.
    \label{lowhigh}
\end{align}
\end{lemma}
\begin{proof}
We establish the first inequality for $\lambda\in(0,1]$ and $r\in[1,+\infty)$, the other cases being handled in the same manner. Indeed, for every $\tau\in I_{j,t}^h$, we have
$j\geq J_\tau+1$, and hence $2^{-j}\leq \var b(\tau)2^{k_0}$. By the definition of the norm and Minkowski's inequality, we have
\begin{align*}
\|f\|_{\widetilde{L}^{r}_{t}(\dot{\B}^{s}_{2,1})}
&=\sum_{j\in\Z} 2^{js}
\left(\int^t_0\|f_j(\tau)\|^r_{L^2}\,d\tau\right)^{1/r}\\
&\lesssim
\sum_{\substack{j\in\Z\\ t_j>0}}
2^{j s}
\left(\int_{I_{j,t}^\ell}
\|f_j(\tau)\|^r_{L^2}\,d\tau\right)^{1/r}
+\sum_{\substack{j\in\Z\\ t_j<t}}
2^{j s}
\left(\int_{I_{j,t}^h}
\|f_j(\tau)\|^r_{L^2}\,d\tau\right)^{1/r}\\
&\lesssim
\|f\|^\ell_{\widetilde{L}^r_{t}(\dot{\B}^{s}_{2,1})}
+\sum_{\substack{j\in\Z\\ t_j<t}}
2^{j (s+1)}
\left(\int_{I_{j,t}^h}
2^{-jr}\|f_j(\tau)\|^r_{L^2}\,d\tau\right)^{1/r}\\
&\leq
\|f\|^\ell_{\widetilde{L}^r_{t}(\dot{\B}^{s}_{2,1})}
+\sum_{\substack{j\in\Z\\ t_j<t}}
2^{j (s+1)}
\left(\int_{I_{j,t}^h}
(\var b(\tau)2^{k_0})^r
\|f_j(\tau)\|^r_{L^2}\,d\tau\right)^{1/r}\\
&\lesssim
\|f\|^\ell_{\widetilde{L}^r_{t}(\dot{\B}^{s}_{2,1})}
+2^{k_0}\var
\sum_{\substack{j\in\Z\\ t_j<t}}
2^{j (s+1)}
\left(\int_{I_{j,t}^h}
b(\tau)^r\|f_j(\tau)\|^r_{L^2}\,d\tau\right)^{1/r}\\
&=
\|f\|^\ell_{\widetilde{L}^r_{t}(\dot{\B}^{s}_{2,1})}
+2^{k_0}
\|\var b(\tau)f\|^h_{\widetilde{L}^r_{t}
(\dot{\B}^{s+1}_{2,1})}.
\end{align*}
\end{proof}

\subsection{A priori estimates}
We start with the following {\emph{a priori}} estimates.
\begin{prop}\label{propapriori}
For any given time $T>0$, let $(n,u)$ be a smooth solution to the Cauchy problem \eqref{eulerre} for $t\in(0,T)$, and let the threshold $J_{t}$ be defined as in \eqref{J}. If 
\begin{align}
    \|(n,\var u)\|_{L^{\infty}_{t}(L^{\infty})}\ll 1,\label{aprioria}
\end{align}
then  $(n,u)$ satisfies
\begin{equation}\label{apriorie}
\mathcal{X} (t)\leq C_{0}{{\big(\theta_\var^{-\frac32}\mathcal{X}_{0}+\theta_\var^{-\frac52}\mathcal{X}^2 (t)\big)}},\quad\quad t\in(0,T),
\end{equation}
where $C_{0}>0$ is a generic constant.
\end{prop}

The proof of Proposition \ref{propapriori} relies on Lemmas \ref{lemmalow} and \ref{lemmahigh} given below. We shall first deal with the underdamped case $0<\lambda\leq 1$, i.e., when $b$ is increasing. The arguments for the overdamped case $\lambda<0$, in which $b$ is decreasing, and for the constant-damping case $\lambda=0$ will be presented at the end of this section.

\subsubsection{Low-frequency analysis for \texorpdfstring{$0< \lambda \leq 1$}{0<lambda<=1}}\label{subsectionlow}

This subsection is devoted to the low-frequency a priori estimates. Following the time-independent case analyzed in \cite{c3}, we introduce a {\emph{generalized damped mode}} in order to partially diagonalize System \eqref{eulerre} in low frequencies. The time-dependent damped mode is defined by
\begin{align}
z:= u+b(t)\nabla n,\label{z}
\end{align}
which can be viewed as a correction associated with Darcy's law, exhibiting stronger regularity and $\mathcal{O}(\varepsilon)$ bounds, which play an essential role in the proof of the relaxation limit. Using this damped mode, we can rewrite \eqref{eulerre} as a heat-type equation with a time-dependent coefficient (consistent with the limiting system \eqref{pm}) coupled with a damped equation, which allows us to establish uniform low-frequency a priori estimates through a pure energy argument.

\begin{lemma}\label{lemmalow}
For any given time $t>0$, let $(n,u)$ be a smooth solution to the Cauchy problem \eqref{eulerre} for $\tau\in(0,t)$ and $0<\lambda\leq 1$, and let the threshold $J_{\tau}$ be defined by \eqref{J}. Then, under the assumption \eqref{aprioria}, it holds
\begin{equation}\label{lowestimate}
    \begin{aligned}
    &\quad \|n\|_{\widetilde{L}^{\infty}_{t}(\dot{\B}^{\frac{d}{2}}_{2,1})}^{\ell}+ \|b(\tau)n\|_{L^1_{t}(\dot{\B}^{\frac{d}{2}+2}_{2,1})}^{\ell}+\var\|u\|_{\widetilde{L}^{\infty}_{t}(\dot{\B}^{\frac{d}{2}}_{2,1})}^{\ell}+\|u\|_{L^1_{t}(\dot{\B}^{\frac{d}{2}+1}_{2,1})}^{\ell}\\
    &\quad\quad  +\|b(\tau)^{-\frac{1}{2}}u\|_{\widetilde{L}^2_{t}(\dot{\B}^{\frac{d}{2}}_{2,1})}^{\ell}+\var^{-1}\|b(\tau)^{-1}z\|_{L^1_{t}(\dot{\B}^{\frac{d}{2}}_{2,1})}^{\ell}\\
    &\lesssim \mathcal{X}_{0}+\mathcal{X}^2(t).
    \end{aligned}
\end{equation}
Here $\mathcal{X} (t)$ and $\mathcal{X}_{0}$ are defined in \eqref{Xp} and \eqref{Xp0}, respectively.   
\end{lemma}

\begin{proof}

With the new unknown $z$, System \eqref{eulerre} can be rewritten as
\begin{equation}\label{eulerlow}
\left\{
    \begin{aligned}
    &\partial_{t}n- P'(\bar{\rho})b(t)\Delta n=-P'(\bar{\rho})\dive z+N_{1},\\
    &\var\partial_{t}z+\frac{1}{\var b(t)}z= \var b'(t)\nabla n+\var b(t)\nabla \partial_{t}n+N_{2},
    \end{aligned}
    \right.
\end{equation}
with the nonlinear terms
$$\left\{
\begin{aligned}
&N_{1}:=-u\cdot\nabla n-G(n)\dive u,\\
&N_{2}:=-\var u\cdot\nabla u.
\end{aligned}
\right.$$
Applying the operator $\dot{\Delta}_{j}$ to $\eqref{eulerlow}_{1}$, taking the scalar product with ${\dot{\Delta}_{j}n}$, integrating it over $\mathbb{R}^{d}$ and employing Bernstein's lemma, we obtain 
\begin{equation*}
\begin{aligned}
&\quad \frac{1}{2}\frac{d}{dt}\|\dot{\Delta}_{j}n\|_{L^{2}}^{2}+\frac{c_{*}2^{2j}}{2} P'(\bar{\rho})b(t)\|\dot{\Delta}_{j}n\|_{L^{2}}^{2}\\
&\leq \Big(\|P'(\bar{\rho})\dive \dot{\Delta}_{j}z\|_{L^2}+\|\dot{\Delta}_{j}N_{1}\|_{L^2} \Big)\|\dot{\Delta}_{j}n\|_{L^2},
\end{aligned}
\end{equation*}
where $c_{*}>0$ is a generic constant. Recall that $I^\ell_{j,t}$ is a connected interval (see Remark \ref{remark11}). By integrating the above inequality over the time interval $I^\ell_{j,t}$ with $j$ satisfying $t_j>0$, it holds that
\begin{equation}\label{lown}
\begin{aligned}
&\quad\sup_{\tau\in I^\ell_{j,t} }\|\dot{\Delta}_{j}n(\tau)\|_{L^2}+\frac{c_{*}2^{2j}}{2}\int_{I^\ell_{j,t}}P'(\bar{\rho})b(\tau)\|\dot{\Delta}_{j}n\|_{L^2}\,d\tau\\
&\leq \|\dot{\Delta}_{j}n_{0}\|_{L^2}+\int_{I^\ell_{j,t}}\Big(\|P'(\bar{\rho})\dive \dot{\Delta}_{j}z\|_{L^2}+\|\dot{\Delta}_{j}N_{1}\|_{L^2}\Big)\,d\tau.
\end{aligned}
\end{equation}
Regarding $z$, one derives from $\eqref{eulerlow}_{2}$ 
that
\begin{equation}\label{lowz1}
    \begin{aligned}
    &\quad\sup_{\tau\in I^\ell_{j,t} }\var\|\dot{\Delta}_{j}z(\tau)\|_{L^2}+\int_{I^\ell_{j,t}}\var^{-1}b(\tau)^{-1}\|\dot{\Delta}_{j}z\|_{L^2}\,d\tau\\
&\leq \var\|\dot{\Delta}_{j}z_{0}\|_{L^2}+\int_{I^\ell_{j,t}}\Big(\var b'(\tau)\|\nabla \dot{\Delta}_{j}n\|_{L^2}+\var b(\tau)\|\partial_\tau\nabla \dot{\Delta}_{j}n\|_{L^2}+\|\dot{\Delta}_{j}N_{2}\|_{L^2}\Big)\,d\tau.
    \end{aligned}
\end{equation}
Using $\eqref{eulerlow}_{1}$ again, we obtain that 
\begin{align*}
  &\quad  \int_{I^\ell_{j,t}}\var b(\tau)\|\partial_\tau\nabla \dot{\Delta}_{j}n\|_{L^2} \,d\tau\\
  &\leqslant \int_{I^\ell_{j,t}}2^{j}\var b(\tau)\Big(2^{2j}b(\tau)P'(\bar{\rho})\|\dot{\Delta}_{j}n\|_{L^2}+P'(\bar{\rho})\|\dive \dot{\Delta}_{j}z\|_{L^2}+\|\dot{\Delta}_{j}N_{1}\|_{L^2}\Big)\,d\tau.
\end{align*}
Since we consider the low frequency regime, we have
\begin{align}
2^{j}\leq (\var b(\tau))^{-1}2^{-k_{0}}.\label{Jsim}
\end{align}
Adding \eqref{lown} to \eqref{lowz1}, we get
\begin{align*}
    &\sup_{\tau\in I_{j,t}^\ell}\Big(\|\dot{\Delta}_{j}n(\tau)\|_{L^2}+\var\|\dot{\Delta}_{j}z(\tau)\|_{L^2}\Big)\\
    &\quad+\frac{c_{*}}{2}P'(\bar{\rho})\int_{I^\ell_{j,t}}2^{2j}b(\tau)\|\dot{\Delta}_{j}n\|_{L^2}\,d\tau+\int_{I^\ell_{j,t}}\var^{-1}b(\tau)^{-1}\|\dot{\Delta}_{j}z\|_{L^2}\,d\tau\\
    \leqslant&\|\dot{\Delta}_{j}n_{0}\|_{L^2}+\var\|\dot{\Delta}_{j}z_{0}\|_{L^2}\\
    &\quad+2^{-k_0}P'(\bar{\rho})\left(\int_{I^\ell_{j,t}}2^{2j}b(\tau)\|\dot{\Delta}_{j}n\|_{L^2}\,d\tau+(2^{-k_0}+1)\int_{I^\ell_{j,t}}\var^{-1}b(\tau)^{-1}\|\dot{\Delta}_{j}z\|_{L^2}\,d\tau\right)\\
    &\quad+\int_{I^\ell_{j,t}} (\var b'(\tau)\|\nabla \dot{\Delta}_{j}n\|_{L^2}+\| \dot{\Delta}_{j}N_1\|_{L^2}+\|\dot{\Delta}_{j}N_{2}\|_{L^2})\,d\tau.
\end{align*}
For the first term on the last line, as we mentioned before, we need to carefully deal with the time derivative of $b(t)$. In fact, we have
\begin{equation}\label{b'1}
\begin{aligned}
    \int_{I^\ell_{j,t}}\var b'(\tau)\|\nabla \dot{\Delta}_{j}n\|_{L^2} \,d\tau
    \lesssim & \sup_{\tau\in I_{j,t}^\ell} \| \dot{\Delta}_{j}n\|_{L^2}\var 2^j\int_{I_{j,t}^\ell} b'(\tau)  \,d\tau\\
    \lesssim & \var b(\min\{t_j,t\})2^j\sup_{\tau\in I_{j,t}^\ell} \| \dot{\Delta}_{j}n\|_{L^2}\\
    \lesssim &2^{-k_0}\sup_{\tau\in I_{j,t}^\ell} \| \dot{\Delta}_{j}n\|_{L^2}.
\end{aligned}
\end{equation}
By choosing $k_0$ large enough, we end up with 
\begin{align*}
    &\sup_{\tau\in I_{j,t}^\ell}(\|\dot{\Delta}_{j}n(\tau)\|_{L^2}+\var\|\dot{\Delta}_{j}z(\tau)\|_{L^2})+\int_{I_{j,t}^\ell}2^{2j}b(\tau)\|\dot{\Delta}_{j}n\|_{L^2}\,d\tau\\
    &\quad +\int_{I_{j,t}^\ell}\var^{-1}b(\tau)^{-1}\|\dot{\Delta}_{j}z\|_{L^2}\,d\tau\\
    \leqslant&C\Big(\|\dot{\Delta}_{j}n_{0}\|_{L^2}+\var\|\dot{\Delta}_{j}z_{0}\|_{L^2}+\int_{I_{j,t}^\ell}(\| \dot{\Delta}_{j}N_1\|_{L^2}+\|\dot{\Delta}_{j}N_{2}\|_{L^2})\,d\tau\Big).
\end{align*}
{
}
Then, multiplying by $2^{j\frac{d}{2}}$ and summing over all $j\in \Z$  satisfying $t_j>0$, we obtain
\begin{equation} \label{lown2}
    \begin{aligned}
    &\quad\|n\|_{\widetilde{L}^{\infty}_{t}(\dot{\B}^{\frac{d}{2}}_{2,1})}^{\ell}+\|b(\tau)n\|_{L^1_{t}(\dot{\B}^{\frac{d}{2}+2}_{2,1})}^{\ell}+\var\|z\|_{\widetilde{L}^{\infty}_{t}(\dot{\B}^{\frac{d}{2}}_{2,1})}^{\ell}+\var^{-1}\|b(\tau)^{-1}z\|_{L^1_{t}(\dot{\B}^{\frac{d}{2}}_{2,1})}^{\ell} \\
    &\lesssim \|(n_{0},\var z_0)\|_{\dot{\B}^{\frac{d}{2}}_{2,1}}+\|(N_{1},N_{2})\|_{L^1_{t}(\dot{\B}^{\frac{d}{2}}_{2,1})}^{\ell},
    \end{aligned}
\end{equation}
   By Lemma \ref{hltransfer}, we know that
   $$
   \var\|z_0\|_{\dot{\B}^{\frac{d}{2}}_{2,1}}\lesssim \var \|u_0\|_{\dot{\B}^{\frac{d}{2}}_{2,1}}+\var \|n_0\|_{\dot{\B}^{\frac{d}{2}+1}_{2,1}}\lesssim \mathcal{X}_0.
   $$

We now estimate the nonlinear terms $N_1$ and $N_2$. It follows from the product law \eqref{uv2} that
    \begin{equation}\label{unablan}
    \begin{aligned}
    \|u\cdot\nabla  n\|_{L^1_{t}(\dot{\B}^{\frac{d}{2}}_{2,1})}\lesssim \|b(\tau)^{-\frac12}u\|_{\widetilde{L}^2_{t}(\dot{\B}^{\frac{d}{2}}_{2,1})}\|b(\tau)^{\frac12}n\|_{\widetilde{L}^2_{t}(\dot{\B}^{\frac{d}{2}+1}_{2,1})}\lesssim \mathcal{X} ^2(t),
          \end{aligned}
\end{equation}  
and
    \begin{equation}\label{unablau}
    \begin{aligned}
    \var\|u\cdot\nabla  u\|_{L^1_{t}(\dot{\B}^{\frac{d}{2}}_{2,1})}\lesssim \|b(\tau)^{-\frac12}u\|_{\widetilde{L}^2_{t}(\dot{\B}^{\frac{d}{2}}_{2,1})}\|\var b(\tau)^{\frac12}u\|_{\widetilde{L}^2_{t}(\dot{\B}^{\frac{d}{2}+1}_{2,1})}\lesssim \mathcal{X} ^2(t),
          \end{aligned}
\end{equation} 
where we have used Lemma \ref{hltransfer} several times, and the $L^2_T$ bound comes from the interpolation between $L^1$ and $L^{\infty}$.
From the product law \eqref{uv2}, the composition estimate \eqref{F1}, and again Lemma \ref{hltransfer}, it also follows that
    \begin{equation}\label{gndivu}
    \begin{aligned}
    \|G(n)\dive u\|_{L^1_{t}(\dot{\B}^{\frac{d}{2}}_{2,1})}&\lesssim \|G(n)\|_{\widetilde{L}^{\infty}_{t}(\dot{\B}^{\frac{d}{2}}_{2,1})}\|u\|_{L^1_{t}(\dot{\B}^{\frac{d}{2}+1}_{2,1})}\\
    &\lesssim \|n\|_{\widetilde{L}^{\infty}_{t}(\dot{\B}^{\frac{d}{2}}_{2,1})}\|u\|_{L^1_{t}(\dot{\B}^{\frac{d}{2}+1}_{2,1})}\\
    &\lesssim \mathcal{X} ^2(t).
          \end{aligned}
\end{equation}  
By \eqref{unablan}-\eqref{gndivu}, we obtain
    \begin{equation}\nonumber
    \begin{aligned}
    \|N_{1}\|_{L^1_{t}(\dot{\B}^{\frac{d}{2}}_{2,1})}^{\ell}&\lesssim\|u\cdot\nabla  n\|_{L^1_{t}(\dot{\B}^{\frac{d}{2}}_{2,1})}^{\ell}+\|G(n)\dive u\|_{L^1_{t}(\dot{\B}^{\frac{d}{2}}_{2,1})}^{\ell}\lesssim\mathcal{X} ^2(t),\\
    \|N_{2}\|_{L^1_{t}(\dot{\B}^{\frac{d}{2}}_{2,1})}^{\ell}&\lesssim \var\|u\cdot\nabla  u\|_{L^1_{t}(\dot{\B}^{\frac{d}{2}}_{2,1})}^{\ell}\lesssim\mathcal{X} ^2(t).
\end{aligned}
\end{equation}  
Substituting the above estimates on $N_{1}$ and $N_{2}$ into \eqref{lown2}, we have
\begin{equation}\label{znlow}
\begin{aligned}
&\quad\|n\|_{\widetilde{L}^{\infty}_{t}(\dot{\B}^{\frac{d}{2}}_{2,1})}^{\ell}+\|b(\tau)n\|_{L^1_{t}(\dot{\B}^{\frac{d}{2}+2}_{2,1})}^{\ell}+\var\|z\|_{\widetilde{L}^{\infty}_{t}(\dot{\B}^{\frac{d}{2}}_{2,1})}^{\ell}+\var^{-1}\|b(\tau)^{-1}z\|_{L^1_{t}(\dot{\B}^{\frac{d}{2}}_{2,1})}^{\ell}\\
 &\lesssim \mathcal{X}_{0}+\mathcal{X} ^2(t).
    \end{aligned}
\end{equation}
In addition, using \eqref{znlow}, $\var\leq1$, and the formula $u=z-b(t)\nabla n$, we recover the following bounds for $u$
\begin{equation*}%\label{ulow}
\begin{aligned}
\var\|u\|_{\widetilde{L}^{\infty}_{t}(\dot{\B}^{\frac{d}{2}}_{2,1})}^{\ell}\lesssim& \var\|z\|_{\widetilde{L}^{\infty}_{t}(\dot{\B}^{\frac{d}{2}}_{2,1})}^{\ell}+\var\|b(\tau)n\|_{\widetilde{L}^{\infty}_{t}(\dot{\B}^{\frac{d}{2}+1}_{2,1})}^{\ell}\\
\lesssim& \var\|z\|_{\widetilde{L}^{\infty}_{t}(\dot{\B}^{\frac{d}{2}}_{2,1})}^{\ell}+2^{-k_0}\|n\|_{\widetilde{L}^{\infty}_{t}(\dot{\B}^{\frac{d}{2}}_{2,1})}^{\ell}\lesssim \mathcal{X}_{0}+\mathcal{X}(t)^2,\\
\|u\|_{L^{1}_{t}(\dot{\B}^{\frac{d}{2}+1}_{2,1})}^{\ell}\lesssim&\var^{-1}  \|b(\tau)^{-1} z\|_{\widetilde{L}^{1}_{t}(\dot{\B}^{\frac{d}{2}}_{2,1})}^{\ell}+\|b(\tau)n\|_{L^{1}_{t}(\dot{\B}^{\frac{d}{2}+2}_{2,1})}^{\ell}\lesssim \mathcal{X}_{0}+\mathcal{X}(t)^2,\\
\|b(\tau)^{-\frac12}u\|_{\widetilde{L}^{2}_{t}(\dot{\B}^{\frac{d}{2}}_{2,1})}^{\ell}\lesssim& \|b(\tau)^{-\frac12}z\|_{\widetilde{L}^{2}_{t}(\dot{\B}^{\frac{d}{2}}_{2,1})}^{\ell}+\|b(\tau)^{\frac12}n\|_{\widetilde{L}^{2}_{t}(\dot{\B}^{\frac{d}{2}+1}_{2,1})}^{\ell}\lesssim \mathcal{X}_{0}+\mathcal{X}(t)^2.
    \end{aligned}
\end{equation*}
These estimates, together with \eqref{znlow}, yield \eqref{lowestimate}.
\end{proof}

\subsubsection{High-frequency analysis  for \texorpdfstring{$0< \lambda \leq 1$}{0<lambda<=1}}\label{subsectionhigh}

Here, we derive the a priori estimates in high frequencies. The proof relies on the construction of localized Lyapunov functionals via hypocoercivity. Compared with the time–independent setting (see e.g. \cite{c2}), there are two essential differences.

First, since the unweighted estimate is not sufficient to control the nonlinear terms, we introduce the time weight $b(t)$ in the $\dot B^{d/2+1}_{2,1}$-estimate. This weighting generates an additional time-dependent linear term involving $b'(t)$, which is absent in the autonomous case and causes a loss of dissipation when $b'(t)>0$. A refined analysis shows that this term can be absorbed by the original dissipation for $0<\lambda<1$ and any $\mu>0$, and for $\lambda=1$ under the condition $\mu>2\varepsilon^2$.

Second, since the norms are restricted to the high--frequency region $j\ge J_t+1$, the time integration must be performed only over $I_{j,t}^h$ (see Remark~\ref{remark11}) for $t>t_j$. Therefore, the value at time $t_j$ necessarily enters the energy estimates,
which requires us to estimate the solution on $[0,t_j]$. Thus, the low and high frequencies are still coupled, and a delicate time decomposition in energy estimates becomes necessary.

\begin{lemma}\label{lemmahigh}
For any given time $t>0$, let $(n,u)$ be any solution to the Cauchy problem \eqref{eulerre} for $\tau\in(0,t)$ and $0<\lambda\leq 1$, and the threshold $J_{t}$ be defined by \eqref{J}. Then, under the assumption \eqref{aprioria}, it holds that
\begin{equation}\label{highestimate}
    \begin{aligned}
&\var\|b(\tau)(n,\var u)\|_{\widetilde{L}^{\infty}_{t}(\dot{\B}^{\frac{d}{2}+1}_{2,1})}^{h}+\frac{\theta_\var}{\var}\|(n,\var u)\|_{L^{1}_{t}(\dot{\B}^{\frac{d}{2}+1}_{2,1})}^{h}+\theta_\var^{\frac{1}{2}}\|b(\tau)^{\frac{1}{2}}(n,\var u)\|_{\widetilde{L}^{2}_{t}(\dot{\B}^{\frac{d}{2}+1}_{2,1})}^{h}\\
&\quad\quad+\theta_\var^{\frac{1}{2}}\|b(\tau)^{-1/2}u\|_{\widetilde{L}^2_{t}(\dot{\B}^{\frac{d}{2}}_{2,1})}^{h}+\frac{\theta_\var}{\var}\|b(\tau)^{-1}u+\nabla n\|_{L^1_{t}(\dot{\B}^{\frac{d}{2}}_{2,1})}^{h}\\
&\quad\lesssim{{ \theta_\var^{-\frac32}\mathcal{X}_{0}+\theta_\var^{-\frac52} \mathcal{X}^2 (t)}},
    \end{aligned}
\end{equation}
where $\mathcal{X} (t)$ and $\mathcal{X}_{0}$ are defined by \eqref{Xp} and \eqref{Xp0}, respectively, and the constant $\theta_\var$ is given in Theorem \ref{theorem11}.
\end{lemma}

\begin{proof}
Applying $\dot{\Delta}_{j}$ to System \eqref{eulerre}, we have
\begin{equation}\label{eulerrej}
\left\{
    \begin{aligned}
    &\partial_{t}\dot{\Delta}_{j}n+u\cdot \nabla \dot{\Delta}_{j} n+P'({\rho})\div \dot{\Delta}_{j} u=R_{1,j},\\
    &\var^2 \partial_{t}\dot{\Delta}_{j} u+\var^2u\cdot \nabla \dot{\Delta}_{j} u+\nabla \dot{\Delta}_{j} n+\frac{1}{b(t)} \dot{\Delta}_{j} u=\var R_{2,j},
    \end{aligned}
    \right.
\end{equation}
with the commutator terms
\begin{equation}\nonumber
\left\{
\begin{aligned}
&R_{1,j}:=u\cdot \nabla \dot{\Delta}_{j} n-\dot{\Delta}_{j}(u\cdot\nabla n)+ G(n)\div \dot{\Delta}_{j} u-\dot{\Delta}_{j}( G(n)\div u),\\
&R_{2,j}:=\var u\cdot \nabla \dot{\Delta}_{j} u-\var \dot{\Delta}_{j}(u\cdot\nabla u).
\end{aligned}
\right.
\end{equation} 
Applying $\nabla$ to $\eqref{eulerrej}_{1}$ and taking the $L^2$-inner with $\dot{\Delta}_{j}\nabla n$, we get
    \begin{equation}\label{ddtnj}
    \begin{aligned}
    &\quad \frac{1}{2}\frac{d}{dt}\| \dot{\Delta}_{j}\nabla n\|_{L^2}^2+\int_{\mathbb{R}^{d}}P'({\rho})\nabla \div \dot{\Delta}_{j} u \cdot \dot{\Delta}_{j}\nabla n\, dx  \\
    &\lesssim  \Big( \| \nabla u \|_{L^{\infty}}\|\dot{\Delta}_{j} \nabla n\|_{L^2}+\| \nabla G(n) \|_{L^{\infty}}\|\div \dot{\Delta}_{j} u\|_{L^2}+\|\nabla R_{1,j}\|_{L^2}\Big)\|\dot{\Delta}_{j}\nabla n\|_{L^2}.
\end{aligned}
\end{equation} 
In order to cancel the second term on the left-hand side of $\eqref{ddtnj}$, applying $\nabla$ to $\eqref{eulerrej}_{2}$, taking the $L^2$-inner with $P'({\rho})\dot{\Delta}_{j} \nabla u$ and then integrating by parts leads to
    \begin{equation}\label{ddtuj}
    \begin{aligned}
    &\quad \frac{\var^2}{2}\frac{d}{dt}\int_{\mathbb{R}^{d}} P'({\rho})|\dot{\Delta}_{j} \nabla u|^2\, dx-\int_{\mathbb{R}^{d}}P'({\rho})\nabla\div \dot{\Delta}_{j} u \cdot\dot{\Delta}_{j}\nabla n\, dx\\
    &\quad+\frac{1}{b(t)}\int_{\mathbb{R}^{d}}  P'({\rho}) |\dot{\Delta}_{j}\nabla u|^2 \, dx\\
    &\lesssim \Big(\|\partial_{t}G(n)\|_{L^{\infty}}+\| \nabla u\|_{L^\infty}+\|\div\big( u P'({\rho}) \big) \|_{L^{\infty}}\Big)\var^2\|\dot{\Delta}_{j}\nabla u\|_{L^2}^2\\
    &\quad+{\var }\|\nabla R_{2,j}\|_{L^2} \|\dot{\Delta}_{j}\nabla u\|_{L^2}+\|\nabla G(n)\|_{L^{\infty}} \|\dot{\Delta}_{j}\nabla u\|_{L^2}\|\dot{\Delta}_{j}\nabla n\|_{L^2}.
    \end{aligned}
\end{equation} 
To capture the dissipation of $\dot{\Delta}_{j}\nabla n$, we shall define a corrector. One derives from  $\eqref{eulerrej}_{1}$-$\eqref{eulerrej}_{2}$ that
    \begin{equation}\label{ddtunj}
    \begin{aligned}
    &\quad \var^2\frac{d}{dt}\int_{\mathbb{R}^{d}} \dot{\Delta}_{j}u\cdot \dot{\Delta}_{j}\nabla n \, dx+\int_{\mathbb{R}^{d}}\big(|\nabla \dot{\Delta}_{j}n|^2-\var^2 P'({\rho})|\div \dot{\Delta}_{j} u|^2 \big) \, dx\\
    &\quad\quad +\int_{\mathbb{R}^{d}}\Big(\frac{1}{b(t)}  \dot{\Delta}_{j} u \cdot \dot{\Delta}_{j} \nabla n+{{\var^2}} u\cdot\nabla \dot{\Delta}_{j}u \cdot \nabla \dot{\Delta}_{j}n -\var^2 u\cdot\nabla \dot{\Delta}_{j}n  \div \dot{\Delta}_{j}u \Big)\,dx\\
    &\leq \var^2\|R_{1,j}\|_{L^2}\|\div \dot{\Delta}_{j}u\|_{L^2}+ \var\|R_{2,j}\|_{L^2} \|\dot{\Delta}_{j} \nabla n\|_{L^2}.
    \end{aligned}
\end{equation}
Now, for a fixed $j$, we shall define the time weighted Lyapunov functional with a parameter $\eta>0$ as
\begin{equation}
\begin{aligned}
    \mathcal{L}_{j,\eta}(t):&=\int_{\mathbb{R}^{d}} \left( b(t)|\dot{\Delta}_{j}\nabla n|^2+b(t)P'({\rho})\var^2|\dot{\Delta}_{j} \nabla u|^2 + \eta \dot{\Delta}_{j}u\cdot \dot{\Delta}_{j}\nabla n \right) \, dx,\label{Lj1}
\end{aligned}
\end{equation}
and the corresponding dissipation rate
\begin{equation}\nonumber
\begin{aligned}
    \mathcal{D}_{j,\eta}(t):&= \var^{-2}\int_{\mathbb{R}^{d}} \left(\eta | \dot{\Delta}_{j}\nabla n|^2+(2-\eta) P'({\rho}) \var^2|\dot{\Delta}_{j} \nabla u|^2+\eta  b(t)^{-1} \dot{\Delta}_{j}u\cdot \dot{\Delta}_{j}\nabla n\right) \, dx.
\end{aligned}
\end{equation}
Choosing $\eta_{1}=2^{-100}\min(2^{-2k_0}P'(\bar{\rho}),1)$, 
with the help of $\|(n,\var u)\|_{L^{\infty}_{t}(L^{\infty})}\lesssim 1$, \eqref{ddtnj}-\eqref{ddtunj} and the fact $2^{-j}\lesssim \var b(t)$, the following Lyapunov inequality holds for $t>t_j$:
\begin{equation}\label{Lyapunov}
\begin{aligned}
&\frac{d}{dt}\mathcal{L}_{j,\eta_1}(t)+\mathcal{D}_{j,\eta_1}(t)-b'(t)\int_{\mathbb{R}^{d}} \left( |\dot{\Delta}_{j}\nabla n|^2+P'({\rho})\var^2|\dot{\Delta}_{j} \nabla u|^2 \right) \, dx\\
\lesssim  &b(t)( \| \nabla u\|_{L^{\infty}}+\|\partial_{t}n\|_{L^{\infty}}+\var^{-1}\| \nabla n\|_{L^{\infty}}) \|\dot{\Delta}_{j}\nabla(n,\var u)\|_{L^2}^2\\
&\quad+b(t)\|\nabla(R_{1,j},R_{2,j})\|_{L^2}\|\dot{\Delta}_{j}\nabla(n, \var u)\|_{L^2}.
\end{aligned}
\end{equation}
To derive a lower bound for the dissipation rate, we first deduce from $\eqref{eulerre}_{1}$ and \eqref{aprioria} that
\begin{align}
 \frac{1}{2}P'(\bar{\rho})\leq P'({\rho})\leq 2P'(\bar{\rho}).\label{Linfty}
\end{align}
Then, by the definition of $t_{j}$, using \eqref{Linfty} and the choice of $\eta_1$, we have, for $t\geqslant t_j$,
\begin{align*}
    &\mathcal{L}_{j,\eta_1}(t)\geqslant \beta b(t)\int_{\mathbb{R}^{d}} \left( |\dot{\Delta}_{j}\nabla n|^2+P'({\rho}) \var^2|\dot{\Delta}_{j} \nabla u|^2 \right) \, dx\sim b(t)\|\nabla \dot{\Delta}_{j}(n,\var u)\|_{L^2}^2,
    \end{align*}
    and
    \begin{align*}
    &\mathcal{D}_{j,\eta_1}(t)\geq \frac{\var^{-2} b(t)^{-1}\eta_1}{2} \mathcal{L}_{j,\eta_1}(t),
    \end{align*}
where, by the choice of $\eta_1$, we may take $\frac12\leq\beta<1$, independently of $j$.

We define $$\widetilde{\sL}_{j,\eta_1}(t):=\sqrt{\sL_{j,\eta_1(t)}/b(t)}$$ such that $\widetilde{\mathcal{L}}_{j,\eta_1}(t)
\sim
\|\nabla\dot{\Delta}_{j}(n,\var u)\|_{L^2}
\sim
2^j\|\dot{\Delta}_{j}(n,\var u)\|_{L^2}$. Combining all the above with Lemma \ref{lem2} leads to 
\begin{align*}
&\sup_{\tau\in[t_j,t]}b(\tau)\widetilde{\mathcal{L}}_{j,\eta_1}(\tau)+\var^{-2}\int_{t_j}^t \widetilde{\mathcal{L}}_{j,\eta_1}(\tau)\,d\tau\\
\lesssim  &\left(\frac{b(t)}{b(t_j)}\right)^{\frac{1+\beta^{-1}}{2}}\bigg(b(t_j)\widetilde{\sL}_{j,\eta_1}(t_j)+\int_{t_j}^t b(\tau)\Big( \| (\nabla u,\partial_{t}n,\var^{-1} \nabla n)\|_{L^{\infty}} \widetilde{\sL}_{j,\eta_1}(\tau)\\
&\quad \quad \quad\quad \quad \quad\quad \quad \quad +\|\nabla(R_{1,j},R_{2,j})\|_{L^2}\Big)\,d\tau\bigg).
\end{align*} 
Since $\lim_{t\rightarrow \infty}b(t)=\infty$, we cannot use this estimate for all $t>t_j$. Hence, we shall cut the time interval again. For a constant $a\in \N$ to be determined later, we set $t_{j,a}>t_j$ such that  $b(t_{j,a})= 2^a b(t_j)$. Then, for $t\in(t_j, t_{j,a})$, we see that 
\begin{align}
&\sup_{\tau\in[t_j,t]}b(\tau)\widetilde{\mathcal{L}}_{j,\eta_1}(\tau)+\var^{-2}\int_{t_j}^t \widetilde{\mathcal{L}}_{j,\eta_1}(\tau)\,d\tau\nonumber\\
\lesssim 2^{\frac{a(1+\beta^{-1})}{2}}  &b(t_j)\widetilde{\sL}_{j,\eta_1}(t_j)+2^{\frac{a(1+\beta^{-1})}{2}}\int_{t_j}^t b(\tau)\left( \| (\nabla u,\partial_{t}n,\var^{-1} \nabla n)\|_{L^{\infty}} \widetilde{\sL}_{j,\eta_1}(\tau)+\|\nabla(R_{1,j},R_{2,j})\|_{L^2}\right)\,d\tau.\nonumber
\end{align} 
For $t>t_{j,a}$, we shall use a new functional. We observe that if we choose $\eta=1$ and then fix an integer $a>k_0+100\max(1,|\log(P'(\bar{\rho}))|)$, we have that, for $t>t_{j,a}$,
\begin{align*}
\int_{\R^d} \dot{\Delta}_j u\cdot \dot{\Delta}_j \nabla n\,dx &\lesssim C 2^{-j}\int_{\R^d}  | \dot{\Delta}_j \nabla u|| \dot{\Delta}_j \nabla n|\,dx\\
= & C\var b(t_j)2^{k_0}\int_{\R^d}  | \dot{\Delta}_j \nabla u|| \dot{\Delta}_j \nabla n|\,dx=C \var b(t_{j,a})2^{k_0-a}\int_{\R^d}  | \dot{\Delta}_j \nabla u|| \dot{\Delta}_j \nabla n|\,dx\\
\lesssim& 2^{k_0-a}b(t)\int_{\R^d}  \var | \dot{\Delta}_j \nabla u|| \dot{\Delta}_j \nabla n|\,dx
\end{align*}
and
\begin{align*}
&\mathcal{L}_{j,1}(t)\geqslant\beta_a b(t)\int_{\mathbb{R}^{d}} \left( |\dot{\Delta}_{j}\nabla n|^2+ P'({\rho})\var^2|\dot{\Delta}_{j} \nabla u|^2 \right) \,dx,\\
    &\mathcal{L}_{j,1}(t)\sim b(t) 2^{2j} \|\dot{\Delta}_{j}(n,\var u)\|_{L^2}^2,\quad\quad \mathcal{D}_{j,1}(t)= \var^{-2} b(t)^{-1} \mathcal{L}_{j,1}(t),
\end{align*}
for some constant $\beta_a>0$ depending only on $a$ and not on $j$. More precisely, the above estimate gives
\begin{align}
\beta_a\geq 1-C2^{k_0-a},\qquad \lim_{a\rightarrow \infty} \beta_a=1. \label{achoosing}
\end{align}
The advantage of choosing $\eta=1$ is that we recover the dissipation rate observed in the spectral analysis.
Then, using the argument above, we end up with 
\begin{align*}&\sup_{\tau\in[t_{j,a},t]}b(\tau)\widetilde{\mathcal{L}}_{j,1}(\tau)+\int_{t_{j,a}}^t (\var^{-2}-2 \beta_a^{-1} b'(\tau))\widetilde{\mathcal{L}}_{j,1}(\tau)\,d\tau\\
\lesssim  &b(t_{j,a})\widetilde{\sL}_{j,1}(t_{j,a})+\int_{t_{j,a}}^t b(\tau)\left(\| (\nabla u,\partial_{t}n,\var^{-1} \nabla n)\|_{L^{\infty}} \widetilde{\sL}_{j,1}(\tau)+\|\nabla(R_{1,j},R_{2,j})\|_{L^2}\right)\,d\tau.
\end{align*}
Now we only need to find an appropriate constant $a$ such that, for $t\geqslant t_{j,a}$, we have 
\begin{align}
\var^{-2}-2\beta_a^{-1}b'(t)>c_*\var^{-2} \theta_{a,\var}>0
\label{condHF}
\end{align}
for some positive constant $c_*<1$, where
\begin{align*}
\theta_{a,\var}:=
\left\{
\begin{array}{ll}
1, & \text{for } 0<\lambda<1,\\[0.3em]
1-\dfrac{2\var^2}{\mu\beta_a}, & \text{for } \lambda=1.
\end{array}
\right.
\end{align*}
Since $\lim_{a\rightarrow \infty} \beta_a=1$, it suffices to have
\begin{align}
\lim_{t\rightarrow\infty}(\var^{-2}-2b'(t))>0.
\label{condHF'}
\end{align}
Here, we need to separate the analysis into two cases.\\
$\bullet$ For $0<\lambda<1$, since
$$
\var^{-2}-2 b'(t)=
\var^{-2}-2\lambda\mu^{-1}(1+t)^{\lambda-1}
$$
increases toward $\var^{-2}$, \eqref{condHF'} is verified for any $\mu>0$. \\
$\bullet$ {{For $\lambda=1$, one observes that
$$
\var^{-2}-2 b'(t)=\var^{-2}-2\mu^{-1},
$$
which corresponds with our assumption $
\mu>2\var^2.$
Without loss of generality, we may assume that $0<\beta_a<1$. In view of
\eqref{achoosing}, we choose $a$ to be the smallest integer satisfying
\begin{align*}
a>k_0+100\max(1,|\log(P'(\bar{\rho}))|),\qquad C2^{k_0-a}\leq\frac{\theta_\var}{2}.
\end{align*}
Then, we have
\begin{align*}
\beta_a\geq \frac12\left(1+\frac{2\var^2}{\mu}\right)
=1-\frac{\theta_\var}{2}.
\end{align*}
Moreover, it follows that
\begin{align*}
\frac12\theta_\var
\leq \theta_{a,\var}
=1-\frac{2\var^2}{\mu\beta_a}
\leq \theta_\var,
\end{align*}
and
\begin{align*}
\var^{-2}-2\beta_a^{-1}b'(t)
=\var^{-2}\left(1-\frac{2\var^2}{\mu\beta_a}\right)
=\var^{-2}\theta_{a,\var}.
\end{align*}
Hence, \eqref{condHF} holds. In addition, since $\frac12\leq\beta<1$ and $a$ is chosen as the smallest integer satisfying the above conditions, we have $2^{\frac{a(1+\beta^{-1})}{2}}\lesssim 2^{\frac{3a}{2}}\lesssim\theta_\var^{-\frac32}$.}}

With the above choice of $a$, we use the functional $\sL_{j,\eta_1}$ for $t\in[t_j,t_{j,a}]$  and $\sL_{j,1}$ for $t\in(t_{j,a},\infty)$.
 
Noticing that $\sL_{j,1}$ and $\sL_{j,\eta_1}$ are equivalent for $t>t_{j,a}$, we end up with, for $t>t_j$,
\begin{align*}
    &\sup_{\tau\in[t_j,t]}b(\tau)\widetilde{\mathcal{L}}_{j,\eta_1}(\tau)+\left\{\begin{array}{ll}\int_{t_j}^t \var^{-2}\widetilde{\mathcal{L}}_{j,\eta_1}\,d\tau &\text{ for } \lambda\in(0,1)\\ \theta_{a,\var}\int_{t_j}^t \var^{-2}\widetilde{\mathcal{L}}_{j,\eta_1}\,d\tau &\text{ for } \lambda=1\end{array}\right.\\
\lesssim  &2^{\frac{a(1+\beta^{-1})}{2}}\bigg(b(t_j)\widetilde{\sL}_{j,\eta_1}(t_j)+\int_{t_j}^t b(\tau)\Big(\| (\nabla u,\partial_{t}n,\var^{-1} \nabla n)\|_{L^{\infty}} \widetilde{\sL}_{j,\eta_1}(\tau)+\|\nabla(R_{1,j},R_{2,j})\|_{L^2}\Big)\,d\tau\bigg).
\end{align*}
Then, multiplying by $\var2^{j\frac d2}$ and summing over all $j\in\Z$ and recalling the above definition of $\theta_{a,\var}$, we have
       \begin{equation}\nonumber
    \begin{aligned}
    &\var\|b(\tau)(n,\var u)\|_{\widetilde{L}^{\infty}_{t}(\dot{\B}^{\frac{d}{2}+1}_{2,1})}^{h}+\var^{-1}\theta_{a,\var}
    \|(n,\var u)\|_{L^{1}_{t}(\dot{\B}^{\frac{d}{2}+1}_{2,1})}^{h}\\
    \lesssim& \theta_\var^{-\frac32}\bigg(\sum_{j\in \Z,t_j<t}\var b(t_j)2^{j\frac d2}\widetilde{\sL}_{j,\eta_1}(t_j)\\ 
    &+\|\var b(\tau) \nabla(n, \var u )\|_{L^{\infty}_{t}(L^{\infty})} \var^{-1}\|(n,\var u)\|_{L^{1}_{t}(\dot{\B}^{\frac{d}{2}+1}_{2,1})}^{h}
    + \| \partial_{t}n\|_{L^1_{t}(L^{\infty})}  \|\var b(\tau)(n,\var u)\|_{\widetilde{L}^{\infty}_{t}(\dot{\B}^{\frac{d}{2}+1}_{2,1})}^{h}\\
    &+\int_{0}^t\sum_{j\geq J_{\tau}-1}2^{(\frac{d}{2}+1)j}\var\|b(\tau) (R_{1,j},R_{2,j})\|_{L^2} \,d\tau\bigg).   \end{aligned}
\end{equation} 
%where $\textbf{1}_{\lambda=1}=1$ if $\lambda=1$ and $0$ otherwise. 

Next, we estimate the nonlinear terms. Due to the embedding $\dot{\B}^{\frac{d}{2}}_{2,1}\hookrightarrow L^{\infty}$, the following estimates hold:
\begin{equation}\nonumber
    \begin{aligned}
    \|\var b(\tau) \nabla(\var   u, n)\|_{L^{\infty}_{t}(L^{\infty})}&\lesssim \var\|b(\tau)  (n,\var u)\|_{L^{\infty}_{t}(\dot{\B}^{\frac{d}{2}+1}_{2,1})}\\
    &\lesssim \|(n,\var u)\|_{\widetilde{L}^{\infty}_{t}(\dot{\B}^{\frac{d}{2}}_{2,1})}^{\ell}+\var\|b(\tau)  (n,\var u)\|_{\widetilde{L}^{\infty}_{t}(\dot{\B}^{\frac{d}{2}+1}_{2,1})}^{h}.
    \end{aligned}
    \end{equation}
Using the equation and the law of products for the low frequencies, we can obtain the bound for $\partial_t n$:
\begin{equation}\nonumber
    \begin{aligned}
   \|\partial_{t}n\|_{L^{1}_t(L^{\infty})}&\lesssim  \|u\|_{L^{1}_t(\dot{\B}^{\frac{d}{2}+1}_{2,1})}+\|u\cdot\nabla n\|_{L^{1}_t(\dot{\B}^{\frac{d}{2}}_{2,1})}+\|G(n)\div u\|_{L^{1}_t(\dot{\B}^{\frac{d}{2}}_{2,1})}\\
   &\lesssim \mathcal{X} (t)+\theta_{a,\var}^{-1}\mathcal{X}^2(t).
    \end{aligned}
    \end{equation} 
For the commutator terms, it holds by Lemma \ref{lemcommutator} that 
\begin{align*}
&\int_{0}^t\sum_{j\geq J_{\tau}-1}2^{(\frac{d}{2}+1)j}\var\|b(\tau) (R_{1,j},R_{2,j})\|_{L^2} \,d\tau\\
\lesssim &\|u\|_{L^1_{t}(\dot{\B}^{\frac{d}{2}+1}_{2,1})}\|\var b(\tau)(n,\var u)\|_{\widetilde{L}^{\infty}_{t}(\dot{\B}^{\frac{d}{2}+1}_{2,1})}+\|b(\tau)^{\frac{1}{2}}n\|_{\widetilde{L}^2_{t}(\dot{\B}^{\frac{d}{2}+1}_{2,1})}\|\var b(\tau)^{\frac{1}{2}}u\|_{\widetilde{L}^{2}_{t}(\dot{\B}^{\frac{d}{2}+1}_{2,1})}\\
    \lesssim& \theta_{a,\var}^{-1}\mathcal{X} ^2(t),
\end{align*}
where we have used the composition estimate 
\begin{equation}\nonumber
\begin{aligned}
\| b(\tau)^{\frac{1}{2}}G(n)\|_{\widetilde{L}^2_{t}(\dot{\B}^{\frac{d}{2}+1}_{2,1})}\lesssim \|b(\tau)^{\frac{1}{2}} n\|_{\widetilde{L}^2_{t}(\dot{\B}^{\frac{d}{2}+1}_{2,1})},
\end{aligned}
\end{equation} 
due to \eqref{aprioria} and Lemma \ref{lemma64}. Then we claim that 
\begin{align*}
\sum_{\substack{j\in\Z\\ j>J_t}}
\var b(t_j)2^{j\frac d2}
\widetilde{\mathcal{L}}_{j,\eta_1}(t_j)
\lesssim
\mathcal{X}_0+\mathcal{X}^2(t).
\end{align*}
We recall that $\widetilde{\mathcal{L}}_{j,\eta_1}(t)
\sim
\|\dot{\Delta}_{j}(\nabla n,\var\nabla u)(t)\|_{L^2}$. For $j\geq J_0+1$, we have $t_j=0$. Thus, it holds
\begin{align*}
\sum_{\substack{j\in\Z\\ j\geq J_0+1}}
\var b(t_j)2^{j\frac d2}
\widetilde{\mathcal{L}}_{j,\eta_1}(t_j)&=
\sum_{\substack{j\in\Z\\ j\geq J_0+1}}
\var b(0)2^{j\frac d2}
\widetilde{\mathcal{L}}_{j,\eta_1}(0)
\\
&\lesssim
\var
\|(n_0,\var u_0)\|_{\dot{\B}^{\frac d2+1}_{2,1}}^{h,J_0}
\lesssim
\mathcal{X}_0.
\end{align*}
For $J_t<j\leq J_0$, by definition, we have $\var b(t_j)2^j=2^{-k_0}$. Then we obtain
\begin{align*}
\sum_{J_t<j\leq J_0}
\var b(t_j)2^{j\frac d2}
\widetilde{\mathcal{L}}_{j,\eta_1}(t_j)
&\lesssim
\sum_{J_t<j\leq J_0}
2^{j\frac d2}
\|\dot{\Delta}_{j}(n,\var u)(t_j)\|_{L^2}
\\
&\lesssim
\|(n,\var u)\|_
{\widetilde{L}^{\infty}_{t}
(\dot{\B}^{\frac d2}_{2,1})}^{\ell}
\lesssim
\mathcal{X}_0+\mathcal{X}^2(t).
\end{align*}
 Here, we used the low-frequency estimate in the last step.

 Gathering all the above estimates, we reach
\begin{equation}\label{highestimate11}
\begin{aligned}
&\var\|b(\tau)(n,\var u)\|_{\widetilde{L}^{\infty}_{t}
(\dot{\B}^{\frac{d}{2}+1}_{2,1})}^{h}+\frac{\theta_{a,\var}}{\var}
\|(n,\var u)\|_{L^{1}_{t}
(\dot{\B}^{\frac{d}{2}+1}_{2,1})}^{h}\lesssim
\theta_\var^{-\frac32}\mathcal{X}_{0}+\theta_\var^{-\frac52}\mathcal{X}^2(t),
\end{aligned}
\end{equation}
from which we infer
\begin{equation}\label{highestimate12}
    \begin{aligned}
  \theta_{a,\var}^{\frac{1}{2}}  \|b(\tau)^{\frac{1}{2}}(n,\var u)\|_{\widetilde{L}^{2}_{t}(\dot{\B}^{\frac{d}{2}+1}_{2,1})}^{h}&\lesssim \Big(\var\|b(\tau)(n,\var u)\|_{\widetilde{L}^{\infty}_{t}(\dot{\B}^{\frac{d}{2}+1}_{2,1})}^{h}\Big)^{\frac{1}{2}}\Big( \frac{ \theta_{a,\var}}{\var}\|(n,\var u)\|_{L^{1}_{t}(\dot{\B}^{\frac{d}{2}+1}_{2,1})}^{h}\Big)^{\frac{1}{2}}\\
    &\lesssim 
 \theta_\var^{-\frac32}\mathcal{X}_{0}+\theta_\var^{-\frac52} \mathcal{X}^2(t).
    \end{aligned}
\end{equation}
And then by \eqref{lowhigh1}-\eqref{lowhigh}, we immediately obtain
\begin{align}
    \theta_{a,\var}^{\frac{1}{2}}  \|b(\tau)^{-\frac{1}{2}} u\|_{\widetilde{L}^{2}_{t}(\dot{\B}^{\frac{d}{2}}_{2,1})}^{h}\lesssim \theta_{a,\var}^{\frac{1}{2}}  \|\var b(\tau)^{\frac{1}{2}} u\|_{\widetilde{L}^{2}_{t}(\dot{\B}^{\frac{d}{2}+1}_{2,1})}^{h}\lesssim  \theta_\var^{-\frac32}\mathcal{X}_{0}+\theta_\var^{-\frac52}\mathcal{X}^2 (t).
\end{align}
Finally, one deduces from \eqref{highestimate11} and the properties in high frequencies that
\begin{equation}\label{highestimate13}
    \begin{aligned}
    \frac{\theta_{a,\var}}{\var}\|b(\tau)^{-1}u+\nabla n\|_{L^1_{t}(\dot{\B}^{\frac{d}{2}}_{2,1})}^{h}\lesssim \theta_{a,\var}\|u\|_{L^1_{t}(\dot{\B}^{\frac{d}{2}+1}_{2,1})}^{h}+\frac{\theta_{a,\var}}{\var}\|n\|_{L^1_{t}(\dot{\B}^{\frac{d}{2}+1}_{2,1})}^{h}\lesssim \theta_\var^{-\frac32}\mathcal{X}_{0}+\theta_\var^{-\frac52}\mathcal{X}^2(t).
    \end{aligned}
\end{equation}
Since $\theta_{a,\var}$ and $\theta_\var$ are uniformly equivalent, we may replace $\theta_{a,\var}$ by $\theta_\var$ on the left-hand sides of \eqref{highestimate11}-\eqref{highestimate13}. Combining \eqref{highestimate11}-\eqref{highestimate13} together, we obtain \eqref{highestimate}. The proof of Lemma  \ref{lemmahigh} is complete.
\end{proof}

\subsubsection{The overdamped case \texorpdfstring{$\lambda<0$}{lambda<0} and the constant damping case \texorpdfstring{$\lambda=0$}{lambda=0}}
Since $b$ is decreasing, for a fixed dyadic block $j$ the low-frequency regime now corresponds to the time-interval $t\in(t_j,\infty)$ and the high-frequency one to $t\in[0,t_j]$. We shall first deal with the high-frequency part. With the same choice of $\eta_1$ as in the previous subsection, we have that for $0<t<t_j$,
\begin{equation}\notag
\begin{aligned}
&\frac{d}{dt}\mathcal{L}_{j,\eta_1}(t)+\mathcal{D}_{j,\eta_1}(t)-b'(t)\int_{\mathbb{R}^{d}} \left( |\dot{\Delta}_{j}\nabla n|^2+ P'({\rho}) \var^2|\dot{\Delta}_{j} \nabla u|^2 \right) \, dx\\
\lesssim  &b(t)\Big( \| \nabla u\|_{L^{\infty}}+\|\partial_{t}n\|_{L^{\infty}}+\var^{-1}\| \nabla n\|_{L^{\infty}}\Big) \|\dot{\Delta}_{j}\nabla(n,\var u)\|_{L^2}^2\\
&\quad +b(t)\|\nabla(R_{1,j},R_{2,j})\|_{L^2}\|\dot{\Delta}_{j}\nabla(n, \var u)\|_{L^2}.
\end{aligned}
\end{equation}
Since $b'(t)<0$, we can directly obtain the desired estimate. For low frequencies, the only change in the analysis is the estimate of the term $\var b'(t)\nabla n$. We have
\begin{align*}
    \int_{t_j}^t\var |b'(\tau)|\|\nabla \dot{\Delta}_{j}n\|_{L^2} \,d\tau
    \leqslant & 4\sup_{\tau\in[t_j,t]} \| \dot{\Delta}_{j}n\|_{L^2}2^j\var\int_{t_j}^t |b'(\tau)| \,d\tau\\
    \leqslant &4 \var b(t_j)2^j\sup_{\tau\in[t_j,t]} \| \dot{\Delta}_{j}n\|_{L^2}\\
    \leqslant&2^{-k_0+2}\sup_{\tau\in[t_j,t]} \| \dot{\Delta}_{j}n\|_{L^2},
\end{align*}
where the second inequality comes from the fact that $b$ is decreasing. The resulting estimate coincides with the definition of our functional space. The overdamped case then follows from the same arguments, with suitable modifications.

The case $\lambda=0$ can be handled by following the underdamped case step by step with straightforward simplifications.

\iffalse
\subsubsection{The final estimate}
    By combining Lemma \ref{lemmalow} with Lemma \ref{lemmahigh}, we obtain \eqref{apriorie}, which concludes the proof of Proposition \ref{propapriori}. Then, {{employing \eqref{Linfty'} and a classical bootstrap argument, one can show that if $\mathcal{X}_0$ is small enough, then  we have \eqref{aprioria}}} and, for all $t\in(0,T)$,
    \begin{align}
    \mathcal{X}(t)\leqslant C \mathcal{X}_0.\label{uniform}
    \end{align}
\fi

\subsubsection{The final estimate}

{{Combining Lemma \ref{lemmalow} with Lemma \ref{lemmahigh}, we obtain
\[
\mathcal X(t)
\leq
C\big(
\theta_\var^{-\frac32}\mathcal X_0
+\theta_\var^{-\frac52}\mathcal X^2(t)
\big).
\]
Assume, as a bootstrap hypothesis, that
\[
\mathcal X(t)\leq 2C\theta_\var^{-\frac32}\mathcal X_0, \qquad t\in(0,T).
\]
Since $\|(n,\varepsilon u)\|_{L^\infty_t(L^\infty)}\lesssim \mathcal{X}(t)$, one can justify \eqref{aprioria} provided $\theta_\var^{-\frac32}\mathcal{X}_0\leq 1$. Then, it holds
\[
\mathcal X(t)
\leq
C\theta_\var^{-\frac32}\mathcal X_0
+
4C^3\theta_\var^{-\frac{11}{2}}\mathcal X_0^2.
\]
If we let $\mathcal X_0$ satisfy the smallness assumption
\[
\mathcal X_0\leq\delta_0\theta_\var^4
\]
for a sufficiently small uniform constant $\delta_0>0$, then we have the uniform bound
\begin{align}\label{uniform}
\mathcal X(t)\leq \frac{3}{2} C\theta_\var^{-\frac32}\mathcal X_0,\qquad t\in(0,T).
\end{align}
}}

\subsection{Proof of global existence and uniqueness}

We first perform the scaling
\[
(\widetilde{n},\widetilde{u})(t,x):= \big(n,\varepsilon u\big)(\var t,x).
\]
The pair $(\widetilde{n},\widetilde{u})$ solves
\begin{equation}\label{eul}
\left\{
\begin{aligned}
&\partial_{t}\widetilde{n}+\widetilde{u}\cdot \nabla \widetilde{n}+\big(P'(\bar{\rho})+G(\widetilde{n})\big)\,\div \widetilde{u}=0,\\
&\partial_{t}\widetilde{u}+\widetilde{u}\cdot \nabla \widetilde{u}+\nabla \widetilde{n}+\frac{1}{\var\, b(\var t)}\, \widetilde{u}=0.
\end{aligned}
\right.
\end{equation}
Since \eqref{eul} is symmetrizable by the matrix 
$$
\begin{pmatrix} (P'(\bar{\rho})+G(\widetilde{n}))^{-1}&0\\[2pt]0& {\rm I}_d\end{pmatrix}
$$
and the damping term $\frac{1}{\var\, b(\var t)}\, \widetilde{u}$ is locally positive in energy estimates, we have the following classical local well-posedness.

\begin{prop}\label{localwell}
For any data $(\widetilde{n}_0,\widetilde{u}_0)\in \B^{\frac{d}{2}+1}_{2,1}$ such that $\inf_{x\in\R^d}\bigl(P'(\bar\rho)+G(\widetilde n_0(x))\bigr)>0$, there exists a time $T_1>0$ such that \eqref{eul} admits a unique classical solution $(\widetilde{n},\widetilde{u})$ with
\[
(\widetilde{n},\widetilde{u})\in \sC^1\big([0,T_1]\times \R^d\big)
\quad\text{and}\quad
(\widetilde{n},\widetilde{u})\in \sC\big([0,T_1];\B_{2,1}^{\frac{d}{2}+1}\big)\cap \sC^1\big([0,T_1];\B_{2,1}^{\frac{d}{2}}\big).
\]
Moreover, if the maximal time of existence $T^*$ is finite, then 
\[
\int_0^{T^*}\|\nabla (\widetilde{n},\widetilde{u})\|_{L^\infty}\, dt=\infty.
\]
\end{prop}

In the following, we denote $W\coloneqq(\widetilde{n},\widetilde{u})$ and $W_0\coloneqq(\widetilde{n}_0,\widetilde{u}_0)$.

\subsubsection*{Step 1: Construction of approximate solutions}
Fix the initial data $W_0\in \dot\B_{2,1}^{\frac{d}{2}}\cap \dot\B_{2,1}^{\frac{d}{2}+1}$ so that the smallness condition from Theorem \ref{theorem11} holds. We approximate the data by
\[
W_0^n \coloneqq ({\rm{Id}}-\dot S_{-n})\,W_0,\qquad n\ge 1.
\]
By construction, $W_0^n\in \B^{\frac{d}{2}+1}_{2,1}$. Consequently, Proposition~\ref{localwell} provides a unique maximal solution
\[
W^n\in \sC\big([0,T^*_n);\,\B_{2,1}^{\frac{d}{2}+1}\big)\cap \sC^1\big([0,T^*_n);\,\B_{2,1}^{\frac{d}{2}}\big).
\]

\subsubsection*{Step 2: Uniform estimates}
Using the {\emph{a priori}} estimate \eqref{uniform} established in the previous section and the fact that
\[
\|W_0^n\|_{\dot\B_{2,1}^{\frac{d}{2}}}+\|W_0^n\|_{\dot\B_{2,1}^{\frac{d}{2}+1}}
\lesssim 
\|W_0\|_{\dot\B_{2,1}^{\frac{d}{2}}}+\|W_0\|_{\dot\B_{2,1}^{\frac{d}{2}+1}},
\]
we infer that
\[
\var\, b(\var \cdot)\, W^n\in L^\infty\big([0,T^*_n);\, \dot\B_{2,1}^{\frac{d}{2}+1}\big).
\]
If $T_n^*<\infty$, then $b(\var\cdot)$ is continuous and positive
on $[0,T_n^*]$. Hence, both $b(\var\cdot)$ and
$b(\var\cdot)^{-1}$ are bounded on this interval. Therefore,
\[
W^n\in L^\infty\big([0,T^*_n);\, \dot\B_{2,1}^{\frac{d}{2}+1}\big)
\quad\text{and thus}\quad
\nabla W^n\in L^\infty\big([0,T^*_n);\, \dot\B_{2,1}^{\frac{d}{2}}\big).
\]
Using the embedding $\dot{\B}_{2,1}^{\frac{d}{2}} \hookrightarrow L^\infty$, the blow-up criterion yields $T^*_n=\infty$.

\subsubsection*{Step 3: Convergence}
\begin{prop}
Let $\widetilde{W}\coloneqq W^1-W^2$, where $W^1$ and $W^2$ are two solutions to \eqref{eul} with initial data $W_0^1$ and $W_0^2$, respectively, and belonging to $\sC(0, T; \dot{\B}_{2,1}^{\frac{d}{2}}\cap\dot{\B}_{2,1}^{\frac{d}{2}+1})$. If both $\|W^1\|_{L^\infty(\dot\B^{\frac{d}{2}}_{2,1})}$ and $\|W^2\|_{L^\infty(\dot\B^{\frac{d}{2}}_{2,1})}$ are smaller than a constant $c>0$, then for all $t\in[0,T]$,
\[
\|\widetilde{W}\|_{L_t^\infty(\dot\B_{2,1}^{\frac{d}{2}})}
\;\lesssim_{c}\;
\|\widetilde{W}_0\|_{\dot\B_{2,1}^{\frac{d}{2}}}
+\int_0^t\!\left(\|(W^1,W^2)\|_{\dot\B_{2,1}^{\frac{d}{2}}}
+\|(W^1,W^2)\|_{\dot\B_{2,1}^{\frac{d}{2}+1}}\right)\,
\|\widetilde{W}\|_{\dot\B_{2,1}^{\frac{d}{2}}}\,d\tau.
\]
\end{prop}

\begin{proof}
The estimate follows from classical stability arguments for symmetric hyperbolic systems in critical Besov spaces, since the damping term does not cause difficulties in the stability estimate.
\end{proof}

Since $W_0^n\to W_0$ in $\dot\B_{2,1}^{\frac{d}{2}}$, the above proposition ensures that $(W^n)_{n\in\N}$ is a Cauchy sequence in $L_T^\infty(\dot\B_{2,1}^{\frac{d}{2}})$ for any finite $T$. Hence, it converges to a limit $W$ in that space. A diagonal extraction then yields
\(
W\in L^\infty(\R_+;\dot\B_{2,1}^{\frac{d}{2}}).
\)
Since the solution is sufficiently regular, passing to the limit in \eqref{eul} is standard. The details are omitted; see, for instance, \cite{c3}.

\subsubsection*{Step 4: Uniqueness}
If $W^1,W^2\in \sC(0, T; \dot{\B}_{2,1}^{\frac{d}{2}}\cap\dot{\B}_{2,1}^{\frac{d}{2}+1})$, then for all $T>0$, using the embedding $L_T^\infty\hookrightarrow L_T^1$ and the continuity of $b$, we have
\[
\int_0^T\!\left(\|(W^1,W^2)\|_{\dot\B_{2,1}^{\frac{d}{2}}}
+\|(W^1,W^2)\|_{\dot\B_{2,1}^{\frac{d}{2}+1}}\right)\,d\tau<\infty.
\]
Moreover, $\|(W^1,W^2)\|_{L^\infty_T(\dot\B_{2,1}^{\frac{d}{2}})}$ is bounded since $W^1,W^2\in \sC(0, T; \dot{\B}_{2,1}^{\frac{d}{2}}\cap\dot{\B}_{2,1}^{\frac{d}{2}+1})$. Combining the stability estimate with Gronwall’s lemma yields $W^1\equiv W^2$ on $[0,T]$. As $T$ is arbitrary, uniqueness holds globally.

Finally, since the pressure law $P$ is smooth, the change of variables between $n^\var$ and $\rho^\var$ is smooth and invertible in a neighborhood of $\bar\rho$. Hence, the estimates for $n^\var$ transfer to $\rho^\var$, which completes the proof of Theorem~\ref{theorem11}.

\section{Strong relaxation limit}\label{sec:relaxation}

\subsection{Derivation of the limit system}

We consider the following approximating compressible Euler system:
\begin{equation}\label{eule}
\left\{
\begin{aligned}
&\partial_{t}\rho^{\var}+\dive (\rho^{\var} u^{\var})=0,\\
&\var^{2} \big( \partial_{t}u^{\var}+u^{\var}\cdot\nabla u^{\var}\big)
+\frac{\nabla P(\rho^{\var})}{\rho^{\var}}
+\frac{u^{\var}}{b(t)}=0,\\
&(\rho^\var,u^\var)(0,x)
=\big(\dot{S}_{J_\var}(\rho^*_0(x)-\bar{\rho})+\bar{\rho},
\var e^{-|x|^2}e_1 \big),
\end{aligned}
\right.
\end{equation}
where $e_1$ denotes the first vector of the canonical basis of $\mathbb R^d$,
and
\[
  J_\var:=J_0=\left[\log_2\frac{\mu}{\var}\right]-k_0
\]
is the initial low-frequency threshold, which tends to $+\infty$ as
$\var\to0$. The initial data satisfy the assumptions of
Theorem~\ref{theorem11} for all sufficiently small $\var$, and the
corresponding global solutions provide an approximating family for the
limiting procedure.

Owing to the uniform bounds obtained in Theorem~\ref{theorem11}, we have
that $\var u^\var$ and $\nabla u^\var$ are uniformly bounded in
$L^\infty(\R^+;\dot{\B}_{2,1}^{\frac d2})$ and
$L^1(\R^+;\dot{\B}_{2,1}^{\frac d2})$, respectively. Hence,
\[
  \var^2 u^\var\cdot\nabla u^\var
  =
  \var(\var u^\var)\cdot\nabla u^\var
  \longrightarrow 0
  \quad\text{in } \mathcal D'(\R^+\times\R^d).
\]
Together with the uniform estimate for the damped mode $z^\var$ defined in
\eqref{z}, this also yields
\[
  \var^2\partial_t u^\var
  \longrightarrow 0
  \quad\text{in } \mathcal D'(\R^+\times\R^d).
\]
Passing to the limit in the second equation of \eqref{eule}, we therefore
obtain
\begin{align*}
    \frac{\nabla P(\rho^{\var})}{\rho^{\var}}
    +\frac{u^{\var}}{b(t)}
    \rightharpoonup 0
    \quad{\rm in}\quad
    \mathcal{D}'(\R^+\times \R^d).
\end{align*}

From the construction of the initial data and the uniform estimates of
Theorem~\ref{theorem11}, we have
\begin{align}
\label{rhoest}
    \|\rho^\var-\bar{\rho}\|_{L^\infty(\R^+;\dot{\B}_{2,1}^{\frac{d}{2}})}
    +
    \|b(\tau)^{\frac12}(\rho^\var-\bar{\rho})\|_{L^2(\R^+;\dot{\B}_{2,1}^{\frac{d}{2}+1})}
    \leqslant
    C\|\rho^*_0-\bar{\rho}\|_{\dot{\B}^{\frac{d}{2}}_{2,1}}
    +C\var .
\end{align}
Moreover,
\[
  \rho^\var(0)-\bar\rho
  =
  \dot S_{J_\var}(\rho_0^*-\bar\rho)
  \longrightarrow
  \rho_0^*-\bar\rho
  \quad\text{in } \dot{\B}^{\frac d2}_{2,1},
\]
since $J_\var\to+\infty$ as $\var\to0$.

In particular, \eqref{rhoest} guarantees the existence of
$\mathcal N$ in
\[
  \bar{\rho}+L^\infty(\R^+;\dot{\B}_{2,1}^{\frac d2})
\]
such that, up to a subsequence,
\begin{align*}
    \rho^\var-\bar{\rho}
    \rightharpoonup
    \mathcal N-\bar{\rho}
    \quad{\rm in}\quad
    \mathcal{D}'(\mathbb{R}^+\times\mathbb{R}^d).
\end{align*}

By the definition of the damped mode, we have
\[
  \rho^\var z^\var
  =
  b(t)\nabla P(\rho^\var)+\rho^\var u^\var .
\]
Inserting this identity into the continuity equation gives
\begin{align*}
    \partial_t\rho^\var-b(t)\Delta P(\rho^\var)=S^\var,
    \qquad
    S^\var:=-\dive(\rho^\var z^\var).
\end{align*}
The uniform estimate on the damped mode further implies that
\[
  S^\var\longrightarrow0
  \quad\text{in }\mathcal D'(\R^+\times\R^d).
\]

On the other hand, one can check that
\[
  \partial_t \rho^\var=-\dive(\rho^\var u^\var)
\]
is uniformly bounded in
$L^1_{\rm loc}(\mathbb{R}^+;\dot{\B}^{\frac d2-1}_{2,1})$.
Therefore, by a compactness argument of Aubin--Lions type, up to extracting
a further subsequence, $\rho^\var-\bar\rho$ converges strongly to
$\mathcal N-\bar\rho$ in $C([0,T];H^{\frac d2-\zeta}_{\rm loc})$ for every $T>0$ and every $\zeta\in(0,1)$. Combining this strong convergence
with the fact that $S^\var\to0$, we may pass to the limit in the equation for
$\rho^\var$ and infer that $\mathcal N=\rho^*$, where $\rho^*$ solves the
porous medium equation
\begin{align}
    \partial_{t}\rho^*-b(t)\Delta P(\rho^*)=0,
    \qquad
    \rho^*(0,x)=\rho_0^*(x).
    \label{eqn}
\end{align}

Furthermore, by lower semicontinuity and \eqref{rhoest}, the solution
$\rho^*$ satisfies
\begin{align}\label{L2:rho*}
 \|\rho^*-\bar{\rho}\|_{L^\infty(\mathbb{R}^+;\dot{\B}_{2,1}^{\frac{d}{2}})}
 +
 \|b(\tau)^{\frac12}(\rho^*-\bar{\rho})\|_{L^2(\mathbb{R}^+;\dot{\B}_{2,1}^{\frac{d}{2}+1})}
 \leq
 C\|\rho^*_0-\bar{\rho}\|_{\dot{\B}^{\frac{d}{2}}_{2,1}}.
\end{align}
Applying the maximal regularity Lemma~\ref{maxiregu} to \eqref{eqn}, we get
\begin{align*}
&\|b(\tau)(\rho^*-\bar{\rho})\|_{L^1(\mathbb{R}^+;\dot{\B}^{\frac{d}{2}+2}_{2,1})}
\\
\lesssim\;&
\|\rho_0^*-\bar{\rho}\|_{\dot{\B}^{\frac{d}{2}}_{2,1}}
+
\|b(\tau)\Delta\big(P(\rho^*)-P(\bar{\rho})
-P'(\bar{\rho})(\rho^*-\bar{\rho})\big)\|_{L^1(\mathbb{R}^+;\dot{\B}^{\frac{d}{2}}_{2,1})}
\\
\lesssim\;&
\|\rho_0^*-\bar{\rho}\|_{\dot{\B}^{\frac{d}{2}}_{2,1}}
+
\|\rho^*-\bar{\rho}\|_{L^\infty(\mathbb{R}^+;\dot{\B}^{\frac{d}{2}}_{2,1})}
\|b(\tau)(\rho^*-\bar{\rho})\|_{L^1(\mathbb{R}^+;\dot{\B}^{\frac{d}{2}+2}_{2,1})}.
\end{align*}
The last term on the right-hand side can be absorbed by the left-hand side,
thanks to \eqref{L2:rho*} and to the smallness of
$\|\rho_0^*-\bar{\rho}\|_{\dot{\B}^{\frac{d}{2}}_{2,1}}$. The $L^1(\mathbb{R}^+;\dot{\B}^{\frac{d}{2}+1}_{2,1})$ bound of $u^*$ is a direct consequence of Darcy's law \eqref{darcy} and the estimates of $\rho^*$.

Finally, uniqueness follows by estimating the difference of two solutions to
\eqref{eqn} with the same initial data in
$L^\infty(\mathbb R^+;\dot{\B}^{\frac d2}_{2,1})$. Since the computations are
standard, we omit the details. Consequently, the whole family converges
toward $\rho^*$, and this completes the proof of
Theorem~\ref{theorem12}.

\begin{remark}
The particular choice of the approximating initial data in \eqref{eule} is
used to ensure that the limit density starts from the prescribed datum
$\rho_0^*$ and that the uniform estimate \eqref{rhoest} depends only on
$\rho_0^*-\bar\rho$, up to a vanishing error as $\var\to0$. For more general families of initial data satisfying the smallness assumptions of
Theorem~\ref{theorem11}, one still obtains uniform bounds analogous to
\eqref{rhoest}. %, with a right-hand side bounded uniformly with respect to $\var$. 
Therefore, the above compactness argument still yields weak
convergence toward a solution of the corresponding
limit equation.%, whose initial datum is determined by the limit of the initial densities.
\end{remark}

However, the above argument only provides weak convergence of the relaxation
limit. To establish global-in-time strong convergence, we need to derive
error estimates between the solutions of the Euler system and the solution
of the limit equation. This is done in the following subsection.

\subsection{Strong convergence to the limit system}
To justify the strong convergence, we first recall that
\begin{align}
    \|\rho^\var-\bar{\rho}\|_{L^\infty_t(\dot{\B}_{2,1}^{\frac{d}{2}})}+\|b(\tau)^{\frac12}(\rho^\var-\bar{\rho})\|_{L^2_t(\dot{\B}_{2,1}^{\frac{d}{2}+1})}+\frac{1}{\var}\|b(\tau)^{-1}z^\var\|_{L^1_t(\dot{\B}_{2,1}^{\frac d2})}\leqslant C_0.\label{uniest}
\end{align}
{{Here, the constant $C_0$ may depend on a positive lower bound for $\theta_\varepsilon$, but we do not track this dependence since $\theta_\varepsilon\to1$ as $\varepsilon\to0$ for any fixed $\mu>0$. }}

We observe that, using similar arguments as in the time-independent case, we can construct a global solution $\rho^*$ of Equation \eqref{eqn}, supplemented with any initial data $\rho^*_0$ such that $\|\rho_0^*-\bar{\rho}\|_{\dot{\B}_{2,1}^{\frac{d}{2}}}$ is small enough. The solution satisfies $\rho^*-\bar{\rho}\in \sC_b(\R_+;\dot{\B}_{2,1}^{\frac{d}{2}})$ and $b(t)(\rho^*-\bar{\rho})\in  L^1(\R_+;\dot{\B}_{2,1}^{\frac{d}{2}+2})$.
Now, we need to separate two distinct cases.

\begin{itemize}
\item {\emph{Overdamped Case $\lambda< 0$.}}
\end{itemize}
We estimate the difference between the solutions of 
\begin{align*}
    \partial_{t}\rho^*-b(t)\Delta P(\rho^*)=0,
\end{align*}
and 
\begin{align*}
    \partial_{t}\rho^{\var}-b(t)\Delta P(\rho^\var)=S^\var.
\end{align*}
 We define $\delta D^\var:=\rho^\var-\rho^*$ that satisfies
\begin{align*}
    \partial_{t}\delta D^\var-b(t)\Delta (P(\rho^\var)-P(\rho^*))=S^\var.
\end{align*}
Using Taylor's formula, there exists a smooth function $\mathcal H^1$ that vanishes at $\bar{\rho}$ such that
\begin{align*}
    P(\rho^\var)-P(\bar{\rho})=P'(\bar{\rho})(\rho^\var-\bar{\rho})+\mathcal H^1(\rho^\var)(\rho^\var-\bar{\rho}),
\end{align*}
and 
\begin{align*}
    P(\rho^*)-P(\bar{\rho})=P'(\bar{\rho})(\rho^*-\bar{\rho})+\mathcal H^1(\rho^*)(\rho^*-\bar{\rho}).
\end{align*}
Now we have 
\begin{align*}
\partial_t\delta D^\var
-
P'(\bar\rho)b(t)\Delta\delta D^\var
&=
b(t)\Delta
\big(
\delta D^\var \mathcal H^1(\rho^\var)
\big)+
b(t)\Delta
\Big(
\big(
\mathcal H^1(\rho^\var)-\mathcal H^1(\rho^*)
\big)
(\rho^*-\bar\rho)
\Big)
+
S^\var.
\end{align*}
Then, using \eqref{maxb}, we obtain
\begin{equation}
\begin{aligned}
   \|\delta D^\var\|_{L^\infty_t(\dot{\B}_{2,1}^{\frac d2-1})}&+\|b(\tau)\delta D^\var\|_{L^1_t(\dot{\B}_{2,1}^{\frac d2+1})}\\
   &\leqslant \|\delta D^\var(0)\|_{\dot{\B}_{2,1}^{\frac d2-1}}+\|S^\var\|_{L^1_t(\dot{\B}_{2,1}^{\frac d2-1})}\\
    &\quad+\|b(\tau)\delta D^\var (\mathcal H^1(\rho^\var)-\mathcal H^1(\bar{\rho}))\|_{L^1_t(\dot{\B}_{2,1}^{\frac d2+1})}\\
    &\quad+\|b(\tau) (\mathcal H^1(\rho^*)-\mathcal H^1(\rho^\var))(\rho^*-\bar{\rho})\|_{L^1_t(\dot{\B}_{2,1}^{\frac d2+1})}.\label{diffest}
\end{aligned}
\end{equation}
To control $S^\varepsilon$, using product laws gives 
\begin{align*}
\|S^\var\|_{L^1_t(\dot{\B}_{2,1}^{\frac d2-1})}\lesssim\|\rho^\var z^\var\|_{L^1_t(\dot{\B}_{2,1}^{\frac d2})}\lesssim\|z^\var\|_{L^1_t(\dot{\B}_{2,1}^{\frac d2})}(\bar{\rho}+\|\rho^\var-\bar{\rho}\|_{L^\infty_t(\dot{\B}_{2,1}^{\frac{d}{2}})}).
\end{align*}
Taking advantage of \eqref{uniest}, we get 
\begin{align}\label{Svar}
    \|S^\var\|_{L^1_t(\dot{\B}_{2,1}^{\frac d2-1})}\lesssim C\var,
\end{align}
where we have used the fact that, for $t>0$,  $b(t)^{-1}\geq b(0)^{-1}>0$ since $\lambda<0$. 

For the two non-linear terms, using the composition and product laws, we have 
\begin{equation}\label{34333}
\begin{aligned}
    &\|b(\tau)\delta D^\var (\mathcal H^1(\rho^\var)-\mathcal H^1(\bar{\rho}))\|_{L^1_t(\dot{\B}_{2,1}^{\frac d2+1})}+\|b(\tau) (\mathcal H^1(\rho^*)-\mathcal H^1(\rho^\var))(\rho^*-\bar{\rho})\|_{L^1_t(\dot{\B}_{2,1}^{\frac d2+1})}\\
    &\lesssim \|b(\tau)\delta D^\var\|_{L^1_t(\dot{\B}_{2,1}^{\frac d2+1})}\|(\rho^\var-\bar{\rho},\rho^*-\bar{\rho})\|_{L^\infty_t(\dot{\B}_{2,1}^{\frac d2})}\\
    &\quad+\|b(\tau)^{\frac 12}\delta D^\var\|_{L^2_t(\dot{\B}_{2,1}^{\frac d2})}\|b(\tau)^{\frac 12}(\rho^\var-\bar{\rho},\rho^*-\bar{\rho})\|_{L^2_t(\dot{\B}_{2,1}^{\frac d2+1})}\\
    &\lesssim \delta_0\Big( \|b(\tau)\delta D^\var\|_{L^1_t(\dot{\B}_{2,1}^{\frac d2+1})}+\|b(\tau)^{\frac 12}\delta D^\var\|_{L^2_t(\dot{\B}_{2,1}^{\frac d2})}\Big),
\end{aligned}
\end{equation}
where we used, thanks to the uniform estimate from Theorem \ref{theorem11}, 
\begin{align*}
    \|(\rho^\var-\bar{\rho},\rho^*-\bar{\rho})\|_{L^\infty_t(\dot{\B}_{2,1}^{\frac d2})}+\|b(\tau)^{\frac 12}(\rho^\var-\bar{\rho},\rho^*-\bar{\rho})\|_{L^2_t(\dot{\B}_{2,1}^{\frac d2+1})}\lesssim \delta_0.
\end{align*}
Inserting \eqref{Svar} and \eqref{34333} into \eqref{diffest} and using 
\begin{align*}
    \|b(\tau)^{\frac12}\delta D^\var\|_{L_t^2
(\dot{\B}_{2,1}^{\frac d2})}
\lesssim
\|\delta D^\var\|_{L_t^\infty
(\dot{\B}_{2,1}^{\frac d2-1})}^{\frac12}
\|b(\tau)\delta D^\var\|_{L_t^1
(\dot{\B}_{2,1}^{\frac d2+1})}^{\frac12},
\end{align*}
we obtain the desired estimate
\begin{align}
    \|\delta D^\var\|_{L^\infty_t(\dot{\B}_{2,1}^{\frac d2-1})}+\|b(\tau)\delta D^\var\|_{L^1_t(\dot{\B}_{2,1}^{\frac d2+1})}\lesssim \|\delta D^\var(0)\|_{\dot{\B}_{2,1}^{\frac d2-1}}+\var.
\end{align}
Concerning $u^\var$, since we have
\begin{align*}
    u^\var-u^*=z^\var-b(t)\Big(\frac{\nabla P(\rho^\var)}{\rho^\var}-\frac{\nabla P(\rho^*)}{\rho^*}\Big),
\end{align*}
the desired bound is obtained by invoking the previously derived error estimate for $\delta D^\var$ and applying the composition law from Lemma~\ref{lemma64}.
\begin{itemize}
\item {\emph{Underdamped Case $0<\lambda\leqslant 1$.}}
\end{itemize}
In this case, we need to estimate the difference in a weighted space. Dividing the equation of $\delta D^\varepsilon$ by $b(t)$ leads to 
\begin{align*}
    \partial_{t}(b(t)^{-1}\delta D^\var)-\Delta (P(\rho^\var)-P(\rho^*))+\frac{b'(t)}{b^2(t)}\delta D^\var=b(t)^{-1}S^\var.
\end{align*}
Because $b'(t)\ge 0$, the contribution $\frac{b'(t)}{b(t)^2}\delta D^\varepsilon$ enters the energy inequality with a favorable sign and may therefore be neglected. We have
\begin{align*}
    &\|b(\tau)^{-1}\delta D^\var\|_{L^\infty_t(\dot{\B}_{2,1}^{\frac d2-1})}+\|\delta D^\var\|_{L^1_t(\dot{\B}_{2,1}^{\frac d2+1})}\\
    &\leqslant \|b(0)^{-1}\delta D^\var(0)\|_{\dot{\B}_{2,1}^{\frac d2-1}}+\|b(\tau)^{-1}S^\var\|_{L^1_t(\dot{\B}_{2,1}^{\frac d2-1})}\\
   &\quad+ \|\delta D^\var (\mathcal H^1(\rho^\var)-\mathcal H^1(\bar{\rho}))\|_{L^1_t(\dot{\B}_{2,1}^{\frac d2+1})}+\| (\mathcal H^1(\rho^*)-\mathcal H^1(\rho^\var))(\rho^*-\bar{\rho})\|_{L^1_t(\dot{\B}_{2,1}^{\frac d2+1})}.\label{diffest1}
\end{align*}
Similarly, we get 
\begin{align*}
\|b(\tau)^{-1}S^\var\|_{L^1_t(\dot{\B}_{2,1}^{\frac d2-1})}\lesssim\|b(\tau)^{-1}\rho^\var z^\var\|_{L^1_t(\dot{\B}_{2,1}^{\frac d2})}\lesssim\|b(\tau)^{-1}z^\var\|_{L^1_t(\dot{\B}_{2,1}^{\frac d2})}(\bar{\rho}+\|\rho^\var-\bar{\rho}\|_{L^\infty_t(\dot{\B}_{2,1}^{\frac{d}{2}})}),
\end{align*}
and
\begin{multline*}
    \|\delta D^\var (\mathcal H^1(\rho^\var)-\mathcal H^1(\bar{\rho}))\|_{L^1_t(\dot{\B}_{2,1}^{\frac d2+1})}+\|(\mathcal H^1(\rho^*)-\mathcal H^1(\rho^\var))(\rho^*-\bar{\rho})\|_{L^1_t(\dot{\B}_{2,1}^{\frac d2+1})}\\
    \lesssim \|\delta D^\var\|_{L^1_t(\dot{\B}_{2,1}^{\frac d2+1})}\|(\rho^\var-\bar{\rho},\rho^*-\bar{\rho})\|_{L^\infty_t(\dot{\B}_{2,1}^{\frac d2})}+\|b(\tau)^{-\frac 12}\delta D^\var\|_{L^2_t(\dot{\B}_{2,1}^{\frac d2})}\|b(\tau)^{\frac 12}(\rho^\var-\bar{\rho},\rho^*-\bar{\rho})\|_{L^2_t(\dot{\B}_{2,1}^{\frac d2+1})}.
\end{multline*}
Then, using the smallness of initial data and interpolation, we end up with 
\begin{align*}
    \|b(\tau)^{-1}\delta D^\var\|_{L^\infty_t(\dot{\B}_{2,1}^{\frac d2-1})}+\|\delta D^\var\|_{L^1_t(\dot{\B}_{2,1}^{\frac d2+1})}\lesssim \|\delta D^\var(0)\|_{\dot{\B}_{2,1}^{\frac d2-1}}+\var,
\end{align*}
and the estimate for $u^\var - u^*$ can be derived in the same way as in the case $\lambda<0$. This establishes the desired error estimates and concludes the proof of Theorem~\ref{theorem13}.     \hfill $\Box$

\section{Appendix}

We begin by recalling the notation associated with the Littlewood–Paley decomposition and Besov spaces. The reader can refer to \cite[Chapter 2]{bahouri1} for a complete overview. We choose a smooth, radial, non-increasing function $\chi(\xi)$ with compact support in $B(0,\frac{4}{3})$ and $\chi(\xi)=1$ in $B(0,\frac{3}{4})$ such that
$$
\varphi(\xi):=\chi(\frac{\xi}{2})-\chi(\xi),\quad \sum_{j\in \mathbb{Z}}\varphi(2^{-j}\xi )=1\,\,\text{for}\,\, \xi\neq 0,\quad \text{{\rm{Supp}}}~ \varphi\subset \{\xi\in \mathbb{R}^{d}~|~\frac{3}{4}\leq |\xi|\leq \frac{8}{3}\}.
$$
For any $j\in \mathbb{Z}$, the homogeneous dyadic blocks $\dot{\Delta}_{j}$ and the low-frequency cut-off operator $\dot{S}_{j}$ are defined by
$$
\dot{\Delta}_{j}u:=\mathcal{F}^{-1}(\varphi(2^{-j}\cdot )\mathcal{F}u),\quad\quad \dot{S}_{j}u:= \mathcal{F}^{-1}( \chi (2^{-j}\cdot) \mathcal{F} u),
$$
where $\mathcal{F}$ and $\mathcal{F}^{-1}$ stand for the Fourier transform and its inverse. Throughout the paper, we may use the notation $\dot{\Delta}_{j}u:=u_{j}.$

Let $\mathcal{S}_{h}'$ be the set of tempered distributions on $\mathbb{R}^{d}$ such that every $u\in \mathcal{S}_{h}'$ satisfies $u\in \mathcal{S}'$ and $\lim_{j\rightarrow-\infty}\|\dot{S}_{j}u\|_{L^{\infty}}=0$. Then, we have
\begin{equation}\nonumber
\begin{aligned}
&u=\sum_{j\in \mathbb{Z}}u_{j}\quad\text{and}\quad \dot{S}_{j}u= \sum_{j'\leq j-1}u_{j'}\quad\text{in}~\mathcal{S}_h'.
\end{aligned}
\end{equation}
With the help of these dyadic blocks, the homogeneous Besov space $\dot{\B}^{s}_{p,r}$, for $p,r\in[1,\infty]$ and $s\in \mathbb{R}$, is defined by
$$
\dot{\B}^{s}_{p,r}:=\{u\in \mathcal{S}_{h}'~|~\|u\|_{\dot{\B}^{s}_{p,r}}:=\|\{2^{js}\|u_{j}\|_{L^p}\}_{j\in\mathbb{Z}}\|_{l^{r}}<\infty\}.
$$

For $T>0, s\in\mathbb{R}$ and $1\leq r,p,\varrho\leq\infty$, we define the homogeneous Chemin-Lerner norm
\begin{align*}
\|u\|_{\widetilde{L}^{\varrho}_{T}(\dot{B}^{s}_{p,r})}:=\big\|(2^{js}\| \dot{\Delta}_{j}u\|_{L^{\varrho}_{T}(L^{p})})\big\|_{l^{r}(\mathbb{Z})}.
\end{align*}
We can then define the homogeneous Chemin-Lerner space $\widetilde{L}_{T}^{\varrho}(\dot{B}^{s}_{p,r})$
as the set of tempered distributions $u$ over $(0,T)\times\mathbb{R}^{d}$ such that $\lim\limits_{j\rightarrow-\infty}\|\dot{S}_{j}u\|_{L^{\varrho}_T(L^{\infty})}=0$ and $\|u\|_{\widetilde{L}_{T}^{\varrho}(\dot{B}^{s}_{p,r})}<\infty$. By the Minkowski inequality, it holds that
\begin{equation}\nonumber
\begin{aligned}
&\|u\|_{\widetilde{L}^{\varrho}_{T}(\dot{\B}^{s}_{p,r})}\leq \|u\|_{L^{\varrho}_{T}(\dot{\B}^{s}_{p,r})}\:\:\text{for}\:\: r\geq\varrho \quad \text{and}\quad \|u\|_{\widetilde{L}^{\varrho}_{T}(\dot{\B}^{s}_{p,r})}\geq \|u\|_{L^{\varrho}_{T}(\dot{\B}^{s}_{p,r})} \:\:\text{for}\:\: r\leq\varrho,
\end{aligned}
\end{equation}
where $\|\cdot\|_{L^{\varrho}_{T}(\dot{\B}^{s}_{p,r})}$ is the usual Lebesgue-Besov norm.

We recall some basic properties of Besov spaces and product estimates which will be used repeatedly in this paper. The reader can refer to \cite[Chapters 2-3]{bahouri1} for more details. The corresponding estimates remain valid in Chemin–Lerner spaces, with the time integrability exponents chosen according to Hölder's
inequality in time.

The first lemma pertains to the so-called Bernstein inequalities.
\begin{lemma}[\!\!\cite{bahouri1}]\label{lemma61}
Let $0<r<R$, $1\leq p\leq q\leq \infty$ and $k\in \mathbb{N}$. For any function $u\in L^p$ and $\lambda>0$, it holds
\begin{equation}\nonumber
\left\{
\begin{aligned}
&{\rm{Supp}}~ \mathcal{F}(u) \subset \{\xi\in\mathbb{R}^{d}~| ~|\xi|\leq \lambda R\}\Rightarrow \|D^{k}u\|_{L^q}\lesssim\lambda^{k+d(\frac{1}{p}-\frac{1}{q})}\|u\|_{L^p},\\
&{\rm{Supp}}~ \mathcal{F}(u) \subset \{\xi\in\mathbb{R}^{d}~|~ \lambda r\leq |\xi|\leq \lambda R\}\Rightarrow \|D^{k}u\|_{L^{p}}\sim\lambda^{k}\|u\|_{L^{p}}.
\end{aligned}
\right.
\end{equation}
\end{lemma}

The following Moser-type product estimates in Besov spaces play a fundamental role in our analysis of nonlinear terms.
\begin{lemma}[\!\!\cite{bahouri1}]\label{lemma63}
The following statements hold:
\begin{itemize}
\item{} Let $p,r\in[1,\infty]$ and $s>0$. Then 
    \begin{equation}\label{uv1}
\begin{aligned}
\|uv\|_{\dot{\B}^{s}_{p,r}}\lesssim \|u\|_{\dot{\B}^{\frac{d}{p}}_{p,1}}\|v\|_{\dot{\B}^{s}_{p,r}}+ \|v\|_{\dot{\B}^{\frac{d}{p}}_{p,1}}\|u\|_{\dot{\B}^{s}_{p,r}}.
\end{aligned}
\end{equation}
\item{}  For $p,r\in[1,\infty]$ and $s\in(-\min\{\frac{d}{p},\frac{d(p-1)}{p}\}, \frac{d}{p}]$, there holds
\begin{equation}\label{uv2}
\begin{aligned}
&\|uv\|_{\dot{\B}^{s}_{p,r}}\lesssim \|u\|_{\dot{\B}^{\frac{d}{p}}_{p,1}}\|v\|_{\dot{\B}^{s}_{p,r}}.
\end{aligned}
\end{equation}
\end{itemize}
\end{lemma}

The following estimates for commutator terms play a role in avoiding the loss of derivatives in high frequencies.

%in $L^{p}$-type Hybrid Besov spaces can be found in \cite{c1}.
\begin{lemma}\label{lemcommutator}
	Let $p,p_1\in [1,\infty ]$  and $p'=\frac{p}{p-1}$. Denote by $[A,B]:=AB-BA$ the commutator bracket.  For $-\min\{\frac{d}{p_1},\frac{d}{p'}\}< s\leq \min\{\frac{d}{p},\frac{d}{p_1}\}+1$, it holds 
        \begin{align}\label{Com:1}
		\sum_{j\in \mathbb{Z}}2^{js}\|[v,\dot{\Delta}_j]\partial_i u\|_{L^p}\lesssim\|\nabla v\|_{\dot{\B}^{\frac{d}{p_1}}_{p_1,1}}\|u\|_{\dot{\B}^{s}_{p,1}}, \ i=1,2,\ldots,d.
	\end{align}
\end{lemma}

We recall a classical estimate regarding the continuity of the composition of functions.
\begin{lemma}\label{lemma64}
Let $d\geq1$, $p,r\in[1,\infty]$, $s>0$ and $F\in C^{\infty}(\mathbb{R})$. Then,  for any $f\in \dot B^s_{p,r}\cap L^\infty$, there exists a constant $C_{f}>0$ depending only on $\|f\|_{L^{\infty}}$, $F$, $s$, $p$ and $d$ such that
\begin{equation}
\begin{aligned}
\|F(f)-F(0)\|_{\dot{\B}^{s}_{p,r}}\leq C_{f}\|f\|_{\dot{\B}^{s}_{p,r}}.\label{F1}
\end{aligned}
\end{equation}
In addition, if $-\min\{\frac dp,\frac{d(p-1)}p\}
<s\leq\frac dp$ and $f_{1}, f_{2}\in\dot{\B}^{s}_{p,r}\cap \dot{\B}^{\frac{d}{p}}_{p,1}$, then we have
\begin{equation}
\begin{aligned}
\|F(f_{1})-F(f_{2})\|_{\dot{\B}^{s}_{p,r}}\leq C_{f_{1},f_{2}}(1+\|(f_{1},f_{2})\|_{\dot{\B}^{\frac{d}{p}}_{p,1}})\|f_{1}-f_{2}\|_{ \dot{\B}^{s}_{p,r}},\label{F3}
\end{aligned}
\end{equation}
where the constant $C_{f_{1},f_{2}}>0$ depends only on $\|(f_{1},f_{2})\|_{L^{\infty}}$, $F$, $s$, $p$ and $d$.
\end{lemma}

We present a lemma that is useful in low-frequency analysis. 
\begin{lemma}
    Let $X:[T_0,T_1]\rightarrow \R_+$ be a continuous function such that $X^2$ is differentiable. Assume that there exist $\sC^1$ functions $c$ and $f$ with $f>0$ and $f'\geq0$ on $[T_0,T_1]$ and a measurable function $A:[T_0,T_1]\rightarrow \R_+$ such that 
    \begin{align*}
        \frac{d}{dt}\bigl(f(t)X(t)^2\bigr)+c(t)X(t)^2
\leqslant A(t)X(t) \quad  a.e.\  on \ [T_0,T_1].
    \end{align*}
    Then, for all $t\in [T_0,T_1]$, we have 
    \begin{align}
     2f(t)X(t)+\int_{T_0}^t(c(\tau)-f'(\tau))X(\tau)\,d\tau\leqslant 2f(T_0)X(T_0)+\int_{T_0}^tA(\tau)\,d\tau,\label{lem2ori}
    \end{align}
    and for any $\alpha>0$ such that $c(\tau)+\alpha f'(\tau)\geq 0$, we have
    \begin{align}
     2f(t)X(t)+\int_{T_0}^t(c(\tau)+\alpha f'(\tau))X(\tau)\,d\tau\leqslant (\frac{f(t)}{f(T_0)})^{\frac{1+\alpha}{2}}\left(2f(T_0)  X(T_0)+\int_{T_0}^tA(\tau)\,d\tau\right).\label{lem2gro}
    \end{align}
    \label{lem2}
\end{lemma}
\begin{proof}
Since the arguments are similar to those of \cite[Lemma 3.1]{danchin5}, we only provide a formal proof here.
Notice that, formally, 
\begin{align*}
  \frac{d}{dt}(f(t)X^2)+cX^2&=f'(t)X^2+2f(t)X\frac{d}{dt}X+cX^2=2X\frac{d}{dt}(fX)+(c-f'(t))X^2.
\end{align*}
Then, formally dividing both sides by $X$ leads to
\begin{align*}
    2\frac{d}{dt}(fX)+(c-f'(t))X\leqslant A.
\end{align*}
Direct integration leads to \eqref{lem2ori}.
For the second one, we rewrite it as 
\begin{align*}
    &\frac{d}{dt}\left(2fX+\int_{T_0}^t(c+\alpha f')X\right)\\
    \leqslant& \frac{1+\alpha}{2}\frac{f'}{f}\left(2fX\right)+A\leqslant \frac{1+\alpha}{2}\frac{f'}{f}\left(2fX+\int_{T_0}^t(c+\alpha f')X\,d\tau\right)+A.
\end{align*}
Then, Gronwall's inequality leads to the desired result \eqref{lem2gro}.
\end{proof}

We then present the time-dependent version of the endpoint maximal regularity. 
\begin{lemma}
\label{maxiregu}
For any given time $T>0$, any nonnegative function $b\in \sC^1([0,T];\R)$ and any $f\in L^1(0,T;\dot{\B}_{2,1}^s)$ with $s\in \R$, let $v$ be a solution to the following Cauchy problem for $t\in(0,T)$:
\begin{equation}
\left\{
    \begin{aligned}
    &\partial_{t}v-b(t)\Delta v=f,\\
    &v(0,x)=v_0(x),
    \label{heatb}
    \end{aligned}
    \right.
\end{equation}
with initial data $v_0\in\dot{\B}_{2,1}^s$. Then, there exists a constant $c_*>0$ such that
\begin{align}
   \| v\|_{\widetilde{L}^\infty_t(\dot{\B}_{2,1}^{s})}+c_*\|b(\tau) v\|_{L^1_t(\dot{\B}_{2,1}^{s+2})}\leqslant \| v_0\|_{\dot{\B}_{2,1}^{s}}+\|f\|_{L^1_t(\dot{\B}_{2,1}^{s})}.\label{maxb}
\end{align}
\end{lemma}
\begin{proof}
Applying the operator $\dot{\Delta}_{j}$ to $\eqref{heatb}_{1}$, taking the scalar product with ${\dot{\Delta}_{j}v}$ and integrating it over $\mathbb{R}^{d}$, we obtain, for $t\in[0,T]$, due to Lemma \ref{lemma61},
\begin{equation}
\begin{aligned}
&\frac{1}{2}\frac{d}{dt}\|\dot{\Delta}_{j}v\|_{L^{2}}^{2}+c_{*}2^{2j} b(t)\|\dot{\Delta}_{j}v\|_{L^{2}}^{2}\leqslant \|\dot{\Delta}_{j}f\|_{L^2}\|\dot{\Delta}_{j}v\|_{L^2},\notag
\end{aligned}
\end{equation}
for some constant $c_{*}=\frac{9}{16}$ related to the support of $\mathcal{F}(\dot{\Delta}_jv)$. Then, it holds by \eqref{lem2ori} that
\begin{equation}
\begin{aligned}
&\|\dot{\Delta}_{j}v(t)\|_{L^2}+ c_{*}2^{2j}\int_{0}^{t}b(\tau)\|\dot{\Delta}_{j}v\|_{L^2}\,d\tau \leqslant \|\dot{\Delta}_{j}v_{0}\|_{L^2}+\int_0^t\|\dot{\Delta}_{j}f\|_{L^2}\,d\tau.\notag
\end{aligned}
\end{equation}
Multiplying by $2^{js}$ and summing over all $j\in\Z$ leads to the desired estimate.
\end{proof}

\vspace{2ex}

\section*{Acknowledgments}

The authors are indebted to the anonymous referees for their valuable comments on the manuscript. T. Crin-Barat is supported by the project ANR-24-CE40-3260 – Hyperbolic Equations, Approximations $\&$ Dynamics (HEAD) and the project ANR-25-CE40-5565 (Cookie).
 X. Pan is supported by the National Natural Science
Foundation of China  under Grant Nos. 12031006 and 12471222. L.-Y. Shou is supported by the National Natural Science
Foundation of China under Grant No. 12301275. Q. Zhu is currently a PhD student, and he would like to thank his supervisor Rapha\"el Danchin for some helpful discussions.

\section*{Data availability statement}

Data sharing is not applicable to this article, as no datasets were generated or analyzed during the current study.

\section*{Conflict of interest statement}

The authors declare that they have no conflict of interest.

\vspace{5mm}

\bigbreak
Timothée Crin-Barat \hfill\break\indent
{\sc Université de Toulouse, Institut de Mathématiques de Toulouse, Route de Narbonne 118, 31062 CEDEX 9 Toulouse, France, \hfill\break\indent
{\it Email address}: {\tt timothee.crin-barat@math.univ-toulouse.fr}}

\bigbreak
Xinghong Pan \hfill\break\indent
{\sc School of Mathematics and Key Laboratory of MIIT, Nanjing University of Aeronautics and Astronautics, Nanjing 211106, China \hfill\break\indent
{\it Email address}: {\tt xinghong\_87@nuaa.edu.cn }}

\bigbreak
Ling-Yun Shou\hfill\break\indent
{\sc School of Mathematical Sciences, Ministry of Education Key Laboratory of NSLSCS, and Key Laboratory of Jiangsu Provincial Universities of FDMTA, Nanjing Normal University, Nanjing 210023, China\hfill\break\indent
{\it Email address}: {\tt shoulingyun11@gmail.com}}

\bigbreak
Qimeng Zhu \hfill\break\indent
{\sc Laboratoire d'Analyse et Math\'ematiques Appliqu\'ees (LAMA UMR8050) 
, \hfill\break\indent
université Paris-Est Cr\'eteil, Cr\'eteil 94010,
France.\hfill\break\indent
{\it Email address}: {\tt qimeng.zhu@u-pec.fr }}

\end{document}